\documentclass[accepted]{uai2025} 
                        
\usepackage[american]{babel}

\usepackage[round]{natbib}

\usepackage{mathtools}
\usepackage{amsmath,amsfonts,amssymb,dsfont,bm,amsthm}
\usepackage{csquotes}
\usepackage{xcolor}
\usepackage{float}
\usepackage{subcaption}
\usepackage{graphicx}
\usepackage{txfonts}
\usepackage{dsfont}
\usepackage{hyperref}
\usepackage{thm-restate} 
\usepackage{soul}
\usepackage{fix-cm}
\usepackage{comment}
\usepackage{soul}

\theoremstyle{definition}

\newtheorem{lemma}{Lemma}

\newtheorem{definition}{Definition}
\newtheorem*{example}{Running Example}

\newtheorem{problem}{Problem}
\newtheorem{assumption}{Assumption}

\usepackage{slashed}

\newcommand{\NC}{\mathcal{NC}}
\newcommand{\CF}{\mathcal{CF}}

\newcommand{\Gf}{\G^f}
\newcommand{\Gs}{\G^s}
\newcommand{\Anc}{\text{Anc}}
\newcommand{\Desc}{\text{Desc}}
\newcommand{\G}{\mathcal{G}}
\newcommand{\C}{\mathcal{C}}
\newcommand{\Do}{\text{do}}
\newcommand{\Forb}{\text{Forb}}
\newcommand{\Ch}{\text{Ch}}
\newcommand{\Xf}{\mathcal{X}^f}
\newcommand{\Xs}{\mathcal{X}^s}

\usepackage{tikz}
\usetikzlibrary{babel}
\usetikzlibrary{arrows}
\usetikzlibrary{decorations}
\usetikzlibrary{decorations.pathmorphing, decorations.pathreplacing, decorations.shapes}

\usepackage[vlined]{algorithm2e}
\RestyleAlgo{ruled}
\SetKwComment{Comment}{/* }{ */}

\definecolor{CentraleRed}{rgb}{0.558,0.09, 0.18}
\definecolor{CentraleGray}{rgb}{0.67,0.67, 0.67}
\definecolor{CentraleBlue}{RGB}{0,67,89}
\definecolor{CentraleGrayEq}{RGB}{66,66,66}

\newcommand{\hiddensection}[1]{
  \addtocontents{toc}{\protect\setcounter{tocdepth}{-1}} 
  \section{#1}
  \addtocontents{toc}{\protect\setcounter{tocdepth}{2}} 
}

\newcommand{\hiddensubsection}[1]{
  \addtocontents{toc}{\protect\setcounter{tocdepth}{-1}} 
  \subsection{#1}
  \addtocontents{toc}{\protect\setcounter{tocdepth}{2}} 
}

\newcommand{\hiddensubsubsection}[1]{
  \addtocontents{toc}{\protect\setcounter{tocdepth}{-1}} 
  \subsubsection{#1}
  \addtocontents{toc}{\protect\setcounter{tocdepth}{2}} 
}

\author[1]{Clément~Yvernes}
\author[1]{Emilie~Devijver}
\author[1]{Eric~Gaussier}

\affil[1]{%
Univ Grenoble Alpes, CNRS, Grenoble INP, LIG
}

\title{Complete Characterization for Adjustment in Summary Causal Graphs of Time Series}

\begin{document}

\maketitle

\begin{abstract}
  The identifiability problem for interventions aims at assessing whether the total causal effect can be written with a do-free formula, and thus be estimated from observational data only. We study this problem, considering multiple interventions, in the context of time series when only an abstraction of the true causal graph, in the form of a summary causal graph, is available. We propose in particular both necessary and sufficient conditions for the adjustment criterion, which we show is complete in this setting, and provide a pseudo-linear algorithm  to decide whether the query is identifiable or not.
\end{abstract}

\hiddensection{Introduction}

Knowing the effect of interventions is key to understanding the effect of a treatment in medicine or the effect of a maintenance operation in IT monitoring systems for example. When one cannot perform interventions in practice, for example when these interventions may endanger people's life or when they may disrupt a critical process or be too costly, one can try and identify do-free formulas which allow one to estimate the effects of interventions using only observational data. 

Finding such do-free formulas is referred to as the identifiability problem for interventions in causal graphs. 
Solving the identifiability problem usually amounts to providing a graphical criterion under which the total effect can be identified, and in providing a do-free formula for its estimation on observational data. 
The problem is, under causal sufficiency, relatively easy for simple graphs,  like DAGs (directed acyclic graphs) for static variables \citep{Pearl_1995} or FTCGs (full time causal graphs) for time series \citep{Blondel_2016}, where the backdoor criterion is sound and complete for monovariate interventions.
It becomes much harder when the graphs considered are abstractions of simple graphs, like CPDAGs (completed partially directed acyclic graphs) and MPDAGs (maximally oriented partially directed) for static variables \citep{Maathuis_2013,Perkovic_2016} or SCGs (summary causal graphs) for time series \citep{assaad_identifiability_2024}. This is due to the fact that, for interventions to be identifiable, one needs to prove that the same do-free formula holds in all the simpler causal graphs corresponding to the abstraction considered.

Despite this increased complexity, \cite{perkovic_identifying_2020} was able to propose, under causal sufficiency, a sound and complete graphical criterion to the identifiability problem for CPDAGs and MPDAGs, namely the general adjustment criterion. However, for SCGs, only a sufficient condition for identifiability has been proposed so far \citep{assaad_identifiability_2024}, using the backdoor criterion which is sound but not complete and under causal sufficiency.

We study in this work the identifiability problem in SCGs under causal sufficiency. In particular:
\begin{itemize}
\item We introduce a common adjustment criterion, which we show is both sound and complete for the adjustment formulae,
\item We propose both necessary and sufficient conditions which further characterize conditions for identifiability by adjustment in SCGs, 
\item Based on these conditions, we derive an algorithm of limited (pseudo-linear) complexity to decide whether the problem is identifiable or not.
\end{itemize}

These results are established here for a single effect and multiple interventions, and hold, on different forms, whether the consistency through time assumption is made or not. {They furthermore rely on novel concepts and tools.}
%

The remainder of the paper is structured as follows: related work is discussed in Section \ref{sec:sota}; Section \ref{sec:notions} introduces the main notions while Section \ref{section:ICA} presents our main result regarding identifiability with the adjustment criterion without assuming \textit{consistency through time}; Section \ref{sec:Consistency_Time} presents a similar result when \textit{consistency through time} holds and a numerical experiment to illustrate the estimation; lastly, Section~\ref{sec:conclusion} concludes the paper. All proofs are provided in the Supplementary Material.

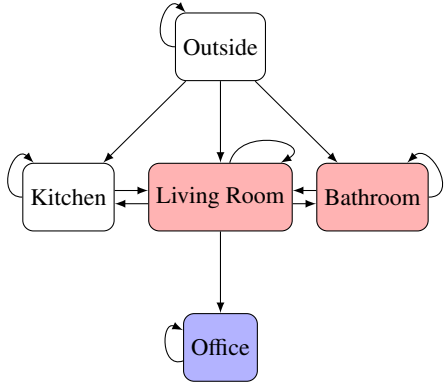
\begin{figure}[t]
  \centering
  \begin{tikzpicture}[{black, rectangle, draw, inner sep=0.1cm}]
    \tikzset{
      nodes={draw, rounded corners, minimum height=0.9cm, minimum width=0.9cm, font=\footnotesize}
    }

    \node (Outside)      at (0,0.5)           {Outside};
    \node[fill=red!30] (LivingRoom) at (0, -1.5) {Living Room};
    \node (Kitchen)      at (-2,-1.5)       {Kitchen};
    \node[fill=red!30] (Bathroom)     at (2,-1.5)        {Bathroom};
    \node[fill=blue!30] (Office)     at (0,-3.5)    {Office};

    \draw[->,>=latex] (Outside)     -- (LivingRoom);
    \draw[->,>=latex] (Outside)     -- (Kitchen);
    \draw[->,>=latex] (Outside)     -- (Bathroom);
    \draw[->,>=latex] (LivingRoom)  -- (Office);

    \begin{scope}[transform canvas={yshift=-.25em}]
      \draw[->,>=latex] (LivingRoom) -- (Kitchen);
    \end{scope}
    \begin{scope}[transform canvas={yshift= .25em}]
      \draw[<-,>=latex] (LivingRoom) -- (Kitchen);
    \end{scope}

    \begin{scope}[transform canvas={yshift=-.25em}]
      \draw[->,>=latex] (LivingRoom) -- (Bathroom);
    \end{scope}
    \begin{scope}[transform canvas={yshift= .25em}]
      \draw[<-,>=latex] (LivingRoom) -- (Bathroom);
    \end{scope}

    \draw[->,>=latex] (Outside)     to [out=180,in=135, looseness=2] (Outside);
    \draw[->,>=latex] (LivingRoom)  to [out= 75,in= 30, looseness=2] (LivingRoom);
    \draw[->,>=latex] (Kitchen)     to [out=180,in=135, looseness=2] (Kitchen);
    \draw[->,>=latex] (Bathroom)    to [out=  0,in= 45, looseness=2] (Bathroom);
    \draw[->,>=latex] (Office)      to [out=200,in=155, looseness=2] (Office);
  \end{tikzpicture}
  \caption{Thermoregulation \citep{assaad_identifiability_2024, Peters_2013}. Only the living room and bathroom have radiators on which we can intervene, highlighted in red. In both scenarios, we are interested in the temperature in the office, highlighted in blue. \textbf{Scenario 1}: \(P(\mathrm{Of}_t \mid \Do(L_{t-1}, L_t, B_{t-1}, B_t))\). \textbf{Scenario 2}: \(P(\mathrm{Of}_t \mid \Do(L_{t-1}, B_{t-1}))\).}
  \label{fig:real_thermoregulation}
\end{figure}

\begin{example}
\label{running_example}
As a running example throughout this paper, we consider the SCG in Figure~\ref{fig:real_thermoregulation}, which models thermoregulation in a house where only the living room and bathroom have radiators. For notational simplicity, let 
\(\,L_t\), \(K_t\), \(\,B_t\), \(\mathrm{Of}_t\), and \(\mathrm{Out}_t\) denote the temperatures in the living room, kitchen, bathroom, office, and outside at time \(t\), respectively. In Scenario~1, we aim to predict the office temperature at time~\(t\), assuming interventions that set the living‐room and bathroom thermostats at times~\(t-1\) and~\(t\): \(P(\mathrm{Of}_t\mid\Do(L_{t-1},L_t,B_{t-1},B_t))\). In Scenario~2, we assume interventions only at time~\(t-1\), giving \(P(\mathrm{Of}_t\mid\Do(L_{t-1},B_{t-1}))\). The Python implementation is available at \href{https://gricad-gitlab.univ-grenoble-alpes.fr/yvernesc/multivariateicainscg}{this repository}.\footnote{https://gricad-gitlab.univ-grenoble-alpes.fr/yvernesc/multivariateicainscg}
\end{example}

\hiddensection{State of the Art}
\label{sec:sota}
%
The identifiability problem for DAGs and under causal sufficiency can be solved with the backdoor criterion, which is sound and complete for total effects with single interventions \citep{Pearl_1995}. However, \cite{Shpitser_2010} have shown that this criterion does not allow one to identify all possible adjustment sets. When the backdoor is not complete, \textit{e.g.}, with hidden confounders or multiple interventions, one may relate to the do-calculus \citep{Pearl_1995} and the associated ID algorithm, which are sound and complete \citep{Shpitser_2010}. 

For CPDAGs, \citet{Maathuis_2013,Perkovic_2016} provided both necessary and sufficient conditions of identifiability  for single interventions, which are nevertheless only sufficient for multiple interventions.
\cite{perkovic_identifying_2020} later developed necessary and sufficient conditions under causal sufficiency and the adjustment criterion for MPDAGs, which encompass DAGs, CPDAGs and CPDAGs with   background knowledge. 
When considering latent confounding, 
\cite{NEURIPS2022_17a9ab41,Wang_2023} provided sufficient conditions of identifiability for PAGs (partial ancestral graphs). 
Cluster DAGs \citep{Anand_2023} constitute another interesting abstraction of simple graphs as they encode partially understood causal relationships between variables grouped into predefined clusters, within which internal causal dependencies remain unspecified. 
They extended do-calculus to establish necessary and sufficient conditions for identifying total effects in these structures.

Fewer studies have however been devoted to the identifiability problem on causal graphs defined over time series, like FTCGs, and abstractions one can define over them \citep{Assaad_2022survey}, like ECGs (extended summary causal graphs) and SCGs. 
As mentioned before, if the problem can be solved relatively easily for FTCGs \citep{Blondel_2016}, it is more complex for SCGs. 
\citet{Eichler_2007} provided sufficient conditions for identifiability of the total effect on graphs based on time series which can directly be generalized to SCGs with no instantaneous relations. With possible instantaneous relations, \cite{Assaad_2023} demonstrated that the total effect is always identifiable on SCGs under causal sufficiency and in the absence of cycles larger than one in the SCG (allowing only self-causes). Another assumption one can make to simplify the problem is to consider that the underlying causal model is linear. This allowed \cite{Ferreira_2024} to propose both necessary and sufficient conditions for identifying direct effects in SCGs. On a slightly different line, \cite{assaad2024identifiabilitytotaleffectssummary} provided sufficient conditions based on the front-door criterion when causal sufficiency is not satisfied. 
The most general result proposed so far on SCGs is the one presented by \cite{assaad_identifiability_2024}, who showed that, under causal sufficiency, the total effect is always identifiable in ECGs and exhibited sufficient conditions for identifiability by common backdoor assuming consistency through time and considering single interventions (but without making assumptions on the form of the SCG or the underlying causal model).

Our work fits within this line of research as it also addresses the identifiability problem in SCGs {and goes further than previous studies by introducing a graphical criterion, shown to be both sound and complete, for identifiability  in SCGs, together with necessary and sufficient conditions allowing one to efficiently decide on identifiability, without other assumptions than causal sufficiency.} 
These results furthermore hold for both single and multiple interventions, with and without consistency through time. 

\hiddensection{Context}
\label{sec:notions}

\hiddensubsection{Notations and Elementary Notions}

For a graph $\G = (\mathcal{V}, \mathcal{E})$, if $X\rightarrow Y$, then $X$ is a \emph{parent} of $Y$ and $Y$ is a \emph{child} of $X$. 
A \emph{path} is a sequence of distinct vertices in which each vertex is connected to its successor by an edge in $\G$.
A \emph{directed path}, or a \textit{causal path}, is a path in which all edges are pointing towards the last vertex. 
A \textit{non-causal} path refers to any path that is not causal. 
If there is a directed path from $X$ to $Y$, then $X$ is an \emph{ancestor} of $Y$, and $Y$ is a \emph{descendant} of $X$. The sets of parents, children, ancestors and descendants of $X$ in $\G$ are denoted by $\text{Pa}(X,\G)$, $\Ch(X,\G)$, $\Anc(X,\G)$ and $\Desc(X,\G)$ respectively. 
We write $X \rightsquigarrow Y$ (or equivalently $Y \leftsquigarrow X$) to indicate that the graph contains a directed path from $X$ to $Y$ consisting of at least one edge.
Furthermore, the \textit{mutilated graph} $\G_{\overline{\textbf{X}} \underline{\textbf{Y}}}$ represents the graph obtained by removing from $\G$ all incoming edge on $\mathbf{X}$ and all outgoing edges from $\mathbf{Y}$.
The \emph{skeleton} of $\G$ is the undirected graph given by forgetting all arrow orientations in $\G$. 
The \textit{subgraph} $\G_{\mid S}$ of a graph $\G$ induced by a vertex set $S$ includes all nodes in $S$ and all edges in $\G$ with both endpoints in $S$.
For two disjoint subsets \textbf{X,Y} \(\subseteq \mathcal{V}\), a \textit{path} from \textbf{X} to \textbf{Y} is a path from some \(X \in \textbf{X}\) to some  \(Y \in \textbf{Y}\). 
A path from \textbf{X} to \textbf{Y} is \textit{proper} if only its first node is in \textbf{X}. 
A \emph{backdoor path} between $X$ and $Y$ is a path between $X$ and $Y$ in which the first arrow is pointing to $X$.
A \emph{directed cycle} is a circular list of distinct vertices in which each vertex is a parent of its successor.  
If a path $\pi$ contains $X_i \rightarrow X_j \leftarrow X_k$ as a subpath, then $X_j$ is a \emph{collider} on $\pi$. A path $\pi$ is \emph{blocked} by a subset of vertices \textbf{Z} if a non-collider in $\pi$ belongs to \textbf{Z} or if $\pi$ contains a collider of which no descendant belongs to \textbf{Z}. Otherwise, \textbf{Z} \emph{d-connects} $\pi$. 

Let $\mathbf{X}, \mathbf{Y}$ and $\mathbf{Z}$ be pairwise distinct sets of variables in a DAG $\G$.
$\mathbf{Z}$ is an \emph{adjustment set} relative to $(\mathbf{X}, \mathbf{Y})$ in $\mathcal{G}$ if for a distribution $P$ compatible with $\mathcal{G}$ 
{\citep[Def. 1.2.2]{Pearl_book2000}} we have\footnote{{As standard in causality studies, $do(\mathbf{x})$ denotes the intervention setting $\mathbf{X}$ to $\mathbf{x}$.}}:
\begin{align*}
    P(\mathbf{y} \mid \Do(\mathbf{x})) = \begin{cases} P(\mathbf{y} \mid \mathbf{x}) & \text{ if } \mathbf{Z} = \emptyset, \\
    \sum_{\mathbf{z}} P(\mathbf{y} \mid \mathbf{x}, \mathbf{z}) P(\mathbf{z}) & \text{ otherwise.} \end{cases}
\end{align*}
Lastly, following \cite{Perkovic_2016}, we make use of the forbidden set in the adjustment criterion.
\begin{definition}[Adjustment criterion]
\label{def:critere_ajustement}
Let $\mathbf{X}, \mathbf{Y}$ and $\mathbf{Z}$ be pairwise distinct sets of variables in a DAG $\G$.
$\mathbf{Z}$ is said to satisfy the \emph{adjustment criterion}  relative to $\mathbf{X}$ and $\mathbf{Y}$  in $\G$ if: 
\begin{enumerate}
    \item \( \Forb\left(\mathbf{X}, \mathbf{Y}, \G\right) \cap \mathbf{Z} = \emptyset\); and
    \item $\mathbf{Z}$ blocks all proper non-causal paths from $\mathbf{X}$  to $\mathbf{Y}$ in $\G$,
\end{enumerate}
where the \emph{forbidden set} \( \Forb\left(\mathbf{X}, \mathbf{Y}, \G\right)\) is the set of all descendants of any \( W \notin \textbf{X}\) which lies on a proper causal path from \textbf{X} to \textbf{Y}.
\end{definition}

\hiddensubsection{Causal Graphs in Time Series}
Consider $\mathcal{V}$ a set of $p$ observational time series and $\mathcal{V}^f=\{\mathcal{V}_{t} | t \in \mathbb{Z}\}$ the set of temporal instances of $\mathcal{V}$ observed over discrete time, where $\mathcal{V}_{t}$ corresponds to the variables of the time series at time $t$. 
We suppose that the discrete time observations $\mathcal{V}^f$ are generated from an \emph{unknown}  structural causal model, which defines an  FTCG which we call the \emph{true} FTCG and a joint distribution $P$ over its vertices which we call the \emph{true} probability distribution, which is compatible with, or Markov relative to, the true FTCG by construction.

As common in causality studies on time series, we consider in the remainder acyclic FTCGs with potential self-causes, \textit{i.e.}, the fact that, for any time series $X$, $X_{t-\ell} \, (\ell \in \mathbb{N}^{*})$ may cause $X_t$. Note that acyclicity is guaranteed for relations between variables at different time stamps and that self-causes are present in most time series. As a result, FTCGs are DAGs in which descendant relationships are constrained by the fact that causality cannot go backward in time, and all causal notions extend directly to FTCGs.

\begin{figure}[t]
\centering

 \begin{subfigure}{0.5\textwidth}
 \centering
 \scalebox{0.80}{
\begin{tikzpicture}[scale=1, transform shape, black, circle, draw, inner sep=0]
 \tikzset{nodes={draw,rounded corners},minimum height=0.6cm,minimum width=0.6cm, font = \tiny}
 \tikzset{latent/.append style={fill=gray!60}}
 
 \node (X) at (1,1) {$X_t$};
 \node (Y) at (1,0) {$Y_t$};
 \node (X-1) at (0,1) {$X_{t-1}$};
 \node (Y-1) at (0,0) {$Y_{t-1}$};
 \node (X-2) at (-1,1) {$X_{t-2}$};
 \node (Y-2) at (-1,0) {$Y_{t-2}$};
 \node (V) at (1,2) {$Z_t$};
 \node (V-1) at (0,2) {$Z_{t-1}$};
 \node (V-2) at (-1,2) {$Z_{t-2}$};
 \draw[->,>=latex] (X-1) to (Y);
 \draw[->,>=latex] (X-2) -- (Y-1);
 \draw[->,>=latex] (X) to  (Y);
 \draw[->,>=latex] (X-1) to (Y-1);
 \draw[->,>=latex] (X-2) to  (Y-2);
 \draw[->,>=latex] (X-2) -- (X-1);
 \draw[->,>=latex] (X-1) -- (X);
 \draw[->,>=latex] (V-2) -- (V-1);
 \draw[->,>=latex] (V-1) -- (V);
 \draw[->,>=latex] (V-2) -- (X-1);
 \draw[->,>=latex] (V-1) -- (X);
 \draw[->,>=latex] (X-2) -- (V-1);
 \draw[->,>=latex] (X-1) -- (V);
 \draw[->,>=latex] (V-2) to  (X-2);
 \draw[->,>=latex] (V-1) to  (X-1);
 \draw[->,>=latex] (V) to  (X);
`
\draw [dashed,>=latex] (V-2) to[left] (-1.55,2);
\draw [dashed,>=latex] (X-2) to[left] (-1.56,1);
 \draw [dashed,>=latex] (Y-2) to[left] (-1.55,0);
 \draw [dashed,>=latex] (V) to[right] (1.55,2);
\draw [dashed,>=latex] (X) to[right] (1.55,1);
\draw [dashed,>=latex] (Y) to[right] (1.55, 0);

\node (a) at (4.65,1) {$X_t$};
\node (b) at (4.65,0) {$Y_t$};
\node (a-1) at (3.65,1) {$X_{t-1}$};
\node (b-1) at (3.65,0) {$Y_{t-1}$};
\node (a-2) at (2.65,1) {$X_{t-2}$};
\node (b-2) at (2.65,0) {$Y_{t-2}$};
\node (c) at (4.65,2) {$Z_t$};
\node (c-1) at (3.65,2) {$Z_{t-1}$};
\node (c-2) at (2.65,2) {$Z_{t-2}$};
 \draw[->,>=latex] (a-1) to (b);
 \draw[->,>=latex] (a-2) -- (b-1);
 \draw[->,>=latex] (a-2) to  (b-2);
 \draw[->,>=latex] (c-2) -- (a-1);
 \draw[->,>=latex] (c-1) -- (a);
 \draw[->,>=latex] (a-1) -- (c);
 \draw[<-,>=latex] (c-1) to  (a-1);
 \draw[->,>=latex] (a-2) -- (a-1);
 \draw[->,>=latex] (a-1) -- (a);
 \draw[->,>=latex] (c-2) -- (c-1);
 \draw[->,>=latex] (c-1) -- (c);
 
\draw [dashed,>=latex] (c-2) to[left] (2.1,2);
\draw [dashed,>=latex] (a-2) to[left] (2.09, 1);
 \draw [dashed,>=latex] (b-2) to[left] (2.1, 0);
\draw [dashed,>=latex] (c) to[right] (5.2,2);
\draw [dashed,>=latex] (a) to[right] (5.2, 1);
\draw [dashed,>=latex] (b) to[right] (5.2, 0);

\node (d) at (8.3,1) {$X_t$};
\node (e) at (8.3,0) {$Y_t$};
\node (d-1) at (7.3,1) {$X_{t-1}$};
\node (e-1) at (7.3,0) {$Y_{t-1}$};
\node (d-2) at (6.3,1) {$X_{t-2}$};
\node (e-2) at (6.3,0) {$Y_{t-2}$};
\node (f) at (8.3,2) {$Z_t$};
\node (f-1) at (7.3,2) {$Z_{t-1}$};
\node (f-2) at (6.3,2) {$Z_{t-2}$};
 \draw[->,>=latex] (d-1) to (e);
 \draw[->,>=latex] (d-2) -- (e-1);
 \draw[->,>=latex] (f-2) -- (d);
 \draw[<-,>=latex] (f-2) to  (d-2);
 \draw[<-,>=latex] (f-1) to  (d-1);
 \draw[<-,>=latex] (f) to  (d);
 \draw[->,>=latex] (d-2) -- (d-1);
 \draw[->,>=latex] (d-1) -- (d);
 \draw[->,>=latex] (f-2) -- (f-1);
 \draw[->,>=latex] (f-1) -- (f);
\draw [dashed,>=latex] (f-2) to[left] (5.75,2);
\draw [dashed,>=latex] (d-2) to[left] (5.74, 1);
 \draw [dashed,>=latex] (e-2) to[left] (5.75, 0);
\draw [dashed,>=latex] (f) to[right] (8.85,2);
\draw [dashed,>=latex] (d) to[right] (8.85, 1);
\draw [dashed,>=latex] (e) to[right] (8.85, 0);

\end{tikzpicture}}
 	\caption{\centering Three FTCGs, $\mathcal{G}_{1}^f$, $\mathcal{G}_{2}^f$ and $\mathcal{G}_{3}^f$.}
 \label{fig:example_FTCG}
 
 \end{subfigure}

  \begin{subfigure}{0.5\textwidth}
 \centering
\begin{tikzpicture}[scale=1, transform shape, black, circle, draw, inner sep=0]
	\tikzset{nodes={draw,rounded corners},minimum height=0.5cm,minimum width=0.5cm, font = \tiny}	
	\tikzset{anomalous/.append style={fill=easyorange}}
	\tikzset{rc/.append style={fill=easyorange}}
 	\node (Z) at (-1.2,0) {$Z$};
	\node (X) at (0,0) {$X$} ;
	\node (Y) at (1.2,0) {$Y$};
 \draw[->,>=latex] (X) -- (Y);
 \begin{scope}[transform canvas={xshift=0, yshift=.1em}]
 \draw [->,>=latex,] (X) -- (Z);
 \end{scope}
 \begin{scope}[transform canvas={xshift=.0, yshift=-.1em}]
 \draw [<-,>=latex,] (X) -- (Z);
 \end{scope}

	\draw[->,>=latex] (X) to [out=155,in=115, looseness=5] (X);
	\draw[->,>=latex] (Z) to [out=155,in=115, looseness=5] (Z);

 \end{tikzpicture}
 \caption{\centering The SCG $\mathcal{G}^s$, reduced from any FTCG in (a).}
 \label{fig:example_SCG}
 \end{subfigure}
 \hfill
 
 \caption{Illustration: (a) three FTCGs; (b) the SCG which can be derived from any FTCG in (a).}
 \label{fig:example_CG}
\end{figure}
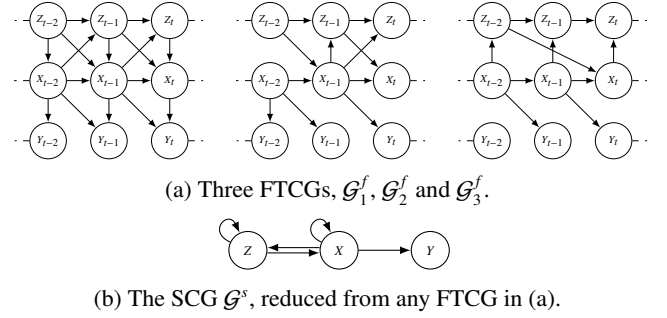%

Experts are used to working with abstractions of causal graphs which summarize the information into a smaller graph that is interpretable, often with the omission of precise temporal information. We consider in this study a known causal abstraction for time series, namely \textit{summary causal graphs}
\citep{Peters_2013, Meng_2020}, which represents causal relationships among time series, regardless of the time delay between the cause and its effect. 
\begin{definition}[Summary causal graph (SCG), Figure~\ref{fig:example_SCG}]
	\label{Summary_G}
	Let $\G^f=(\mathcal{V}^f,\mathcal{E}^f)$ be an FTCG built from the set of time series $\mathcal{V}$. The \emph{summary causal graph} (SCG) $\G^s=(\mathcal{V}^s, \mathcal{E}^s)$ associated to $\G^f$ is such that:
     \begin{itemize}
         \item $\mathcal{V}^s$ corresponds to the set of time series $\mathcal{V}$,
         \item $X \rightarrow Y \in \mathcal{E}^s$ if and only if there exists at least one timepoint $t$ and one temporal lag $0\leq \gamma$ such that $X_{t-\gamma} \rightarrow Y_t \in \mathcal{E}^f$.
     \end{itemize} 
In that case, we say that $\G^s$ is \emph{reduced from} $\G^f$. 
\end{definition}
SCGs may include directed cycles and even self-loops. For example, the three FTCGs in Figure \ref{fig:example_FTCG} are acyclic, while the SCG in Figure \ref{fig:example_SCG} has a cycle. 
We use the notation $X \rightleftarrows Y$ to indicate situations where there exist time instants in which $X$ causes $Y$ and $Y$ causes $X$. It is furthermore worth noting that if there is a single SCG reduced from a given FTCG, different FTCGs, with possibly different orientations and skeletons, can yield the same SCG. 
For example, the SCG in Figure \ref{fig:example_SCG}  can be reduced from any FTCG in Figure \ref{fig:example_FTCG}, even though they may have different skeletons or  different orientations.
In the remainder, we refer to any FTCG from which a given SCG $\Gs$ can be reduced as a \emph{candidate FTCG} for $\Gs$.
 For example, in Figure \ref{fig:example_CG}, $\G_{1}^f$, $\G_{2}^f$ and $\G_{3}^f$ are all candidate FTCGs for $\G^s$. The class of all candidate FTCGs for $\Gs$ is denoted by $\C(\Gs)$.

\hiddensubsection{Problem Setup}
We focus in this paper on identifying total effects \citep{Pearl_book2000} of multiple interventions on {single effects, written $P \left(Y_{t} = y_{t} \middle\vert \bigl(\Do\bigl(X^i_{t_i} = x^i_{t_i}\bigr) \bigr)_i \right)$ (as well as $P\left(y_{t} \mid \Do \left((x^i_{t_i})_i\right)\right)$ by a slight abuse of notation)} when only the SCG reduced from the true FTCG is known, using the common adjustment criterion defined below.

\begin{definition}[Common adjustment criterion]
\label{def:IBC}
    Let $\Gs = (\mathcal{V}^s, \mathcal{E}^s)$ be an SCG. 
    Let $\mathbf{X}, \mathbf{Y}$ and $\mathbf{Z}$ be pairwise distinct subsets of $\mathcal{V}^f$. 
    $\mathbf{Z}$ satisfies the \emph{common adjustment criterion} relative to $\mathbf{X}$ and $\mathbf{Y}$ in $\Gs$ if for all $\Gf \in \mathcal{C}(\Gs)$, 
     $\mathbf{Z}$ satisfies the adjustment criterion relative to $\mathbf{X}$ and $\mathbf{Y}$ in $\Gf$.
\end{definition}
%
%
This criterion is sound and complete for the adjustment formulae, meaning that:
\begin{restatable}{proposition}{mypropadjustmentSCG}
\label{prop:adjustment_scg}
    Let $\Gs = (\mathcal{V}^s, \mathcal{E}^s)$ be an SCG and let $\mathbf{X}, \mathbf{Y}$ and $\mathbf{Z}$ be pairwise distinct subsets of $\mathcal{V}^f$. We say that a probability distribution $P$ is \emph{compatible} with $\mathcal{G}^s$ if there exists $\mathcal{G}^f \in \C(\Gs)$ such that $P$ is compatible with $\Gf$. The two following propositions are equivalent:
    \begin{description}
        \item[(i)] $\mathbf{Z}$ satisfies the common adjustment criterion relative to $\mathbf{X}$ and $\mathbf{Y}$,
        \item[(ii)]  for all $P$ compatible with $\Gs$:
    \begin{equation}
    \label{eq:adjustment}
         P\left(\mathbf{y} \mid \Do(\mathbf{x})\right)
        = \begin{cases}
            P\left(\mathbf{y} \mid \mathbf{x}\right) & \text{if } \mathbf{Z} = \emptyset\\
            \sum_{\mathbf{z}} P\left(\mathbf{y} \mid \mathbf{x}, \mathbf{z}\right) P(\mathbf{z}) & \text{otherwise.}
            \end{cases}
    \end{equation}
\end{description}
When either (i) or (ii) hold, we say that the total effect $P\left(\mathbf{y} \mid \Do(\mathbf{x}) \right)$ is \emph{identifiable in $\Gs$ by adjustment criterion}.
\end{restatable}

{Finally, our problem takes the form:}

\begin{problem}
\label{pb:mains_objective}
    Consider an SCG $\Gs$. 
    We aim to find out {operational}\footnote{{That is, conditions one can rely on in practice. In particular the number of candidate FTCGs in $\C(\Gs)$ is usually too costly to enumerate (it may even be infinite); operational conditions should thus not rely on the enumeration of all FTCGs.}} necessary and sufficient conditions to identify the total effect $P\left(y_{t} \mid \Do \left((x^i_{t_i})_i\right)\right)$ by common adjustment when having access solely to the SCG $\Gs$.

\end{problem}
Note that if $Y$ is not a descendant of one of the intervening variables $X^i$ in $\Gs$ {or if $\gamma_i \coloneqq t - t_i < 0$}, then $X^i_{t_i}$ can be removed from the conditioning set through, \textit{e.g.}, the adjustment for direct causes \citep{Pearl_book2000}.  In the extreme case where $Y$ is not a descendant of any element of $\{X^i\}_i$, then $P \left(y_t \mid \Do( x^1_{t-\gamma_1}), \dots , \Do( x^n_{t-\gamma_n}) \right) = P(y_t)$. In the remainder, we thus assume that $Y$ is a descendant of each element in $\{X^i\}_i$ in $\Gs$ {and that $\gamma_i \geq 0$ for all $i$}, and will use the following notations: $\Xf \coloneqq \{X^i_{t-\gamma_i}\}_i$ and $\Xs \coloneqq \{ X^i\}_i$. 

\hiddensection{Identifiability by Common Adjustment}\label{section:ICA}

 We provide in this section the main results of this paper, which is a graphical necessary and sufficient condition  for identifiability of the causal effect by common adjustment, and a solution to compute it in practice. The classical \textit{consistency through time}, assuming that causal relations are the same at different time instants, is not assumed here and its discussion is postponed to Section \ref{sec:Consistency_Time}. All the proofs are deferred to Section \ref{sec:proof:4} in the Supplementary Material.
 
\hiddensubsection{Necessary and Sufficient Condition Based on the Common Forbidden Set}\label{subsection:equiv_based_CD}

We first introduce the \emph{common forbidden set}, the set of vertices that belong to \( \Forb\left(\Xf, Y_t, \G^f\right)\) in at least one candidate FTCG \( \G^f \). The common forbidden set, and the related notion of non-conditionable set defined below, define a set of variables which cannot be elements of a common adjustment set as they violate the first condition in    Definition~\ref{def:critere_ajustement}. As such, they cannot be used as conditioning variables in the do-free formula rewriting the interventions (Equation~\ref{eq:adjustment} in Proposition~\ref{prop:adjustment_scg}).

\begin{definition}
\label{def:CD}
   Let $\Gs = (\mathcal{V}^s, \mathcal{E}^s)$ be an SCG and $P(y_t \mid \text{do}(x^1_{t-\gamma_1}), \dots, \text{do}(x^n_{t-\gamma_n}))$ be the considered effect. We define the \emph{common forbidden set} as follows:
    $$\CF \coloneqq \bigcup_{\Gf \in \mathcal{C}(\Gs)} \Forb\left(\Xf, Y_t, \Gf\right).$$
The \emph{set of non-conditionable variables}  is defined by 
    $$\NC \coloneqq \CF \setminus \Xf.$$
\end{definition}

\begin{example}
\label{ex:3}
In the first scenario, we have $\NC = \left\{ \mathrm{Of}_{t- 1}, \mathrm{Of}_t \right\}$, whereas, in the second scenario, we have $\NC = \left\{K_{t-1}, K_t, B_t, L_t, \mathrm{Of}_{t-1}, \mathrm{Of}_t \right\}$. In the second scenario, $K_{t-1}$ cannot belong to a common adjustment set as there exists a candidate FTCG which contains the path $L_{t-1} \rightarrow K_{t-1} \rightarrow L_t \rightarrow \mathrm{Of}_t$. Similarly, $L_t$ cannot belong to a common adjustment set as there exists a candidate FTCG which contains the path $L_{t-1} \rightarrow L_t\rightarrow  \mathrm{Of}_t$.
\end{example}

Theorem~\ref{th:equiv_IBC_multivarie} below shows that identifiability by common adjustment is directly related to the existence of collider-free backdoor path remaining in this set. 

\begin{restatable}{theorem}{mytheoremequivIBC}{}
    \label{th:equiv_IBC_multivarie}
    Let $\Gs = (\mathcal{V}^s, \mathcal{E}^s)$ be an SCG and $P(y_t \mid \Do(x^1_{t-\gamma_1}), \dots, \Do(x^n_{t-\gamma_n}))$ be the considered effect. Then the two statements are equivalent:
    \begin{enumerate}
        \item The effect is identifiable by common adjustment in $\Gs$. \label{th:equiv_IBC_multivarie:1}
        \item For all intervention $X^i_{t-\gamma_i}$ and candidate FTCG $\Gf \in \C(\Gs)$, $\Gf$ does not contain a collider-free backdoor path going from $X^i_{t-\gamma_i}$ to $Y_t$ that remains in $\NC \cup \{ X^i_{t-\gamma_i}\}$. \label{th:equiv_IBC_multivarie:2}
    \end{enumerate}
    
    In that case, a common adjustment set is given by $\C \coloneqq \left(\mathcal{V}^f \setminus \NC \right) \setminus X^f$, and we have 
    \begin{equation*}
        P(y_t \mid \Do((x^i_{t-\gamma_i})_i))
        = \sum_{\mathbf{c}} P\left(y_t  \mid (x^i_{t-\gamma_i})_i, \mathbf{c}\right) P(\mathbf{c}).
    \end{equation*}
\end{restatable}

\begin{figure}[t]
    \centering
    \begin{tikzpicture}[scale = 1.5]
        \coordinate (Y) at (0,0);
        \coordinate (Xi) at (-1.5,3);
        \coordinate (Xj) at (.8,1.25);
        \coordinate (C) at (.35,1.5);
        \coordinate (D) at (.35,1);
        
        \draw[dashed] (Xi) -- (-2,3);
        \draw[dashed] (-2,3) -- (-2,0);
        \draw[dashed] (Xi) -- (-.2,3);
        \draw[dashed] (-.2,3) -- (-.2,1.25);
        
        \draw[dashed] (Xj) -- (-.2,1.25);
        \draw[dashed] (Xj) -- (1.3,1.25);
        \draw[dashed] (1.3,1.25) -- (1.3,0);

        \draw[->] (-2.5,3.25) -- (-2.5,-.25);
        \fill (-2.5,3)  circle (.5pt) node[left] {\small $t-\gamma_i$};
        \fill (-2.5,1.25)  circle (.5pt) node[left] {\small $t-\gamma_j$};
        \fill (-2.5,0)  circle (.5pt) node[left] {\small $t$};
        
        \draw[green!50!black, thick, decorate, decoration={snake, pre length = 2 mm, amplitude=.4mm, segment length=3mm, post=lineto, post length=1mm}, -] (Xi) .. controls (-1.5,2) and (.435,2.5) .. (C) ;
        \draw[green!50!black, thick, {latex[scale=.9]}-] (D) -- (C);
        \draw[green!50!black, thick, decorate, decoration={snake, amplitude=.4mm, segment length=3mm, post length = 1 mm}] (D) .. controls (.35,.7) and (0,.3) .. (Y);

        \fill (Y)  circle (1pt) node[below] {$Y_t$};
        \fill (Xi) circle (1pt) node[above] {$X^i_{t-\gamma_i}$};
        \fill (Xj) circle (1pt) node[above right] {$X^j_{t-\gamma_j}$};
        \fill (C)  circle (1pt) node[above right] {\small $C_{t_c}$};
        \fill (D)  circle (1pt) node[below right] {\small $D_{t_d}$};
    \end{tikzpicture}
    \caption{Proof idea of Theorem~\ref{th:equiv_IBC_multivarie}. The green path represents a proper non-causal path $\pi^f$ from $\Xf$ to $Y_t$ starting at $X^i_{t-\gamma_i}$. The node $X^j_{t-\gamma_j}$ represents another intervention (if any). The dashed lines depict the set $\CF \cup \{ X^i_{t-\gamma_i} \}$. The vertex $C_{t_c}$ is the last node on $\pi^f$ outside of $\CF \cup \{ X^i_{t-\gamma_i} \}$, and $D_{t_d}$ is its successor on $\pi^f$. Necessarily, $\pi^f$ must contain the arrow $C_{t_c} \rightarrow D_{t_d}$; otherwise, $C_{t_c}$ would belong to $\CF$.}\label{fig:CD_proof}
\end{figure}
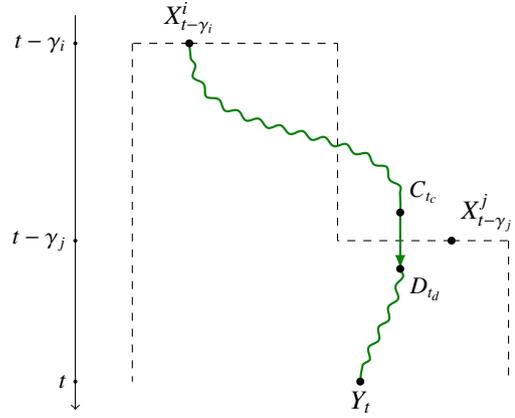

\begin{proof}[Proof Sketch]
Let $\Gf$ be a candidate FTCG. For any $X^i_{t-\gamma_i} \in \Xf$, consider any proper non-causal path $\pi^f$ from $\Xf$ to $Y_t$ that starts at $X^i_{t-\gamma_i}$. Then $\pi^f$ either:
\begin{itemize}
    \item leaves $\CF \cup \{X^i_{t-\gamma_i}\}$, in which case it contains a non-collider in $\C$ (see $C_{t_c}$ in Figure \ref{fig:CD_proof}) and is blocked by $\C$,
    \item remains in $\CF \cup \{X^i_{t-\gamma_i}\}$ and contains a collider, in which case it is also blocked by $\C$, since the collider and its descendants remain in $\NC$,
    \item remains in $\CF \cup \{X^i_{t-\gamma_i}\}$ and contains no collider (i.e., it is a collider-free backdoor path), in which case it cannot be blocked.
\end{itemize}
Thus, collider-free backdoor paths entirely contained in $\NC \cup \{X^i_{t-\gamma_i}\}$ are therefore the only proper non-causal paths from $\Xf$ to $Y_t$ starting at $X^i_{t-\gamma_i}$ that cannot be blocked by $\C$. Moreover, such paths cannot be blocked by any common adjustment set, as they remain in $\NC \cup \Xf$. As a result, the effect is identifiable by common adjustment if and only if no such path exists for any $X^i_{t-\gamma_i} \in \Xf$.
\end{proof}

\begin{example}
In the first scenario, we have \(\NC = \left\{ \mathrm{Of}_{t-1}, \mathrm{Of}_t \right\}\). No candidate FTCG contains a non-causal path that remains within \(\NC\). As a result, \(P(\mathrm{Of}_t \mid \Do(L_{t-1}, L_t, B_{t-1}, B_t))\) is identifiable by common adjustment. In the second scenario,we have $\NC = \left\{K_{t-1}, K_t, B_t, L_t, \mathrm{Of}_{t-1}, \mathrm{Of}_t \right\}$ and we know that both \(K_{t-1}\) and \(L_t\) cannot be part of a common adjustment set. Since a candidate FTCG contains the path \(L_{t-1} \leftarrow K_{t-1} \rightarrow L_t \rightarrow \mathrm{Of}_t\), \(P(\mathrm{Of}_t \mid \Do(L_{t-1}, B_{t-1}))\) is not identifiable by common adjustment.
\end{example}

\hiddensubsection{An Efficient Way to Decide {on Identifiability}} \label{ssct:in practice}

To determine whether the causal effect is identifiable, we propose an algorithm that efficiently tests the existence of collider-free backdoor paths that remain in $\NC$, except perhaps for their first vertices. Instead of enumerating all FTCGs and such paths in $\NC$, which would be computationally prohibitive, we introduce a more refined approach to characterize their existence. Specifically, we distinguish between those with and without forks. Paths without forks can be easily and efficiently identified. The situation is more complex for those that contain forks, but they can still be efficiently handled via a divide-and-conquer strategy. All these elements are detailed in the following subsections.

\subsubsection{Additional Characterizations}

 \paragraph{Characterization of \texorpdfstring{$\NC$}{NC}}

We first introduce another characterization of $\NC$ based on the time instant a time series first arrives in this set.

\begin{definition}
    Let $\Gs = (\mathcal{V}^s, \mathcal{E}^s)$ be an SCG and $P(y_t \mid \text{do}(x^1_{t-\gamma_1}), \dots, \text{do}(x^n_{t-\gamma_n}))$ be the considered effect. For a time series $F \in \mathcal{V}^s$, we define
    \[
    t_{\NC}(F) \coloneqq \min \{ t_1 \mid F_{t_1} \in \NC \},
    \]
    as the first time step at which $F$ enters the non-conditionable set $\NC$, with the convention that $\min \{\emptyset\} = +\infty$.
\end{definition}

\begin{example}
In the second scenario,  we have $\NC = \left\{K_{t-1}, K_t, B_t, L_t, \mathrm{Of}_{t-1}, \mathrm{Of}_t \right\}$. As a result, $t_{\NC}(\mathrm{Kitchen}) = t-1$ and $t_{\NC}(\mathrm{Outside}) = +\infty$.
\end{example}

In the next lemma, we show that $\{t_{\NC}(F)\}_{F \in \mathcal{V}^S}$ gives a simple characterization of these sets.

\begin{restatable}{lemma}{mylemmadefequivCD}{(Characterization of $\NC$)}
\label{lemma:def_equiv_CD}
Let $\Gs = (\mathcal{V}^s, \mathcal{E}^s)$ be an SCG and let $P(y_t \mid \Do(x^1_{t-\gamma_1}), \dots, \text{do}(x^n_{t-\gamma_n}))$ be the considered effect. With the convention $\{F_{t_1}\}_{t_1 \geq +\infty} = \emptyset$, we have:
\[
 \NC =  \bigcup_{Z \in \mathcal{V}^S} \{Z_{t_1}\}_{t_1 \geq t_{\NC}(Z)} \setminus \Xf.
\]
Moreover, $(t_{\NC}(F))_{F \in \mathcal{V}^s}$ {can be} computed through Algorithm~\ref{algo:calcul_t_NC}, {detailed} in Appendix \ref{sec:proof:4}, {which} complexity is 
pseudo-linear with respect to $\Gs$ and $\Xf$.
\end{restatable}
{The above characterization, based on $t_{\NC}(F)$, slightly departs from standard, purely graphical characterizations often used in the identifiability literature. This is due to the complexity of the class of candidate FTCGs and the difficulty to explore this class efficiently.}

\begin{algorithm}[t]
\setlength{\algomargin}{0.5em} 
\setlength{\rightskip}{-.5cm} 
\caption{Computation of $(t_{\NC}(S))_{S \in \mathcal{V}^s}$}\label{algo:calcul_t_NC}
\SetKwInput{KwData}{Input}
\SetKwInput{KwResult}{Output}
\KwData{ $\Gs = (\mathcal{V}^s, \mathcal{E}^s)$ an SCG and $\Xf$.}
\KwResult{$(t_{\NC}(S))_{S \in \mathcal{V}^s}$}

\tcp{Compute $t_C \quad \forall ~ C \in \Ch(\Xs)$. (cf. Lemma \ref{lemma:characterisation_CF})} 

$AncY \gets \left(\max \left\{ t_1 \mid  \exists\Gf\text{ s.t } 
  S_{t_1} \in \Anc(Y_{t}, \Gf \setminus \Xf)  \right\}\right)_{S \in \mathcal{V}^s}$ \;
\ForEach{$C \in \Ch(\Xs)$}{
    $t_{\min} \gets \min \{t-\gamma_i \mid X^i \in \text{Pa}(C,\Gs)\}$\;
    $t_C \gets \min \left\{ t_1 \in \left[ t_{\min}, AncY[C] \right]  \mid C_{t_1} \notin \Xf  \right\}$\;
    $d(C) \begin{aligned}[t]
        \gets \# \{i \mid X^i \in \text{Pa}(C, \Gs) \text{ and } t - \gamma_i < t_C \} \geq 1& \\
        \textbf{or } \# \{i \mid X^i \in \text{Pa}(C, \Gs) \text{ and } t - \gamma_i = t_C \} \geq 2 &\text{ ;}
    \end{aligned}$ 
}

\tcp{Avoid extra computations. (cf. Lemma \ref{lemma:characterisation_CF_suite})}

$L \gets [(C,t_C)]_{C \in \Desc(\Xf), \text{with } t_C < +\infty}$ \;
Sort $L$ using $(t_C, not~ d(C))$ lexicographically\;

\tcp{Compute $(t_{\NC}(S))_{S \in \mathcal{V}^s}$. (cf. Lemma \ref{lemma:preuve_algo_t_NC})}

$(t_{\NC}(S)) \gets +\infty \quad \forall S \in \mathcal{V}^s$\;
$S.seen \gets False \quad \forall S \in \mathcal{V}^s$\;
\For{$(C,t_C) \in L$}{
    \eIf{d(C)}{
        \ForEach{unseen $D \in \Desc(C,\Gs)$}{
            $t_{\NC}(D) \gets \min \{t_1 \mid t_1 \geq t_C \text{ and } D_{t_1} \notin \Xf\}$\;
            $D.seen \gets \text{true}$ \;
        }
    }{
        \ForEach{unseen $D \in \Desc(C,\Gs\setminus \Xs)$}{
            $t_{\NC}(D) \gets \min \{t_1 \mid t_1 \geq t_C \text{ and } D_{t_1} \notin \Xf\}$\;
            $D.seen \gets \text{true}$ \;
        }
        \ForEach{unseen $D \in \Desc(C,\Gs)$}{
            $t_{\NC}(D) \gets \min \{t_1 \mid t_1 \geq t_C+1 \text{ and } D_{t_1} \notin \Xf\}$\;
            $D.seen \gets \text{true}$ \;
        }
    }
}
\end{algorithm}

\paragraph{Collider-free backdoor paths without fork in \texorpdfstring{$\NC$}{NC}}
\label{subsubsection:CfbWF}
First, Lemma \ref{lemma:IBC_enumeration_chemins_diriges}  characterizes efficiently the existence of collider-free backdoor paths that do not contain a fork.

\begin{restatable}{lemma}{mylemmaIBCenumerationcheminsdiriges}{(Characterization of collider-free backdoor paths without fork)}
\label{lemma:IBC_enumeration_chemins_diriges}
   Let $\Gs = (\mathcal{V}^s, \mathcal{E}^s)$ be an SCG and $P(y_t \mid \text{do}(x^1_{t-\gamma_1}), \dots, \text{do}(x^n_{t-\gamma_n}))$ be the considered effect. The following statements are equivalent:
    \begin{enumerate}
        \item There exists an intervention $X^i_{t-\gamma_i}$ and a candidate FTCG $\Gf \in \C(\Gs)$ which contains  $X^i_{t- \gamma_i} \leftsquigarrow Y_t$ which remains in $\NC \cup \{X^i_{t- \gamma_i}\}$.
        \label{lemma:IBC_enumeration_chemins_diriges:1}
        \item There exists an intervention $X^i_{t-\gamma_i}$ such that $ \gamma_i = 0$ and $X^i \in \Desc \left(Y, \Gs_{\mid \mathcal{S}} \right)$, where $\mathcal{S} \coloneqq \{S \in \mathcal{V}^s \mid t_{\NC}(S) \leq t \}\cup \{X^i \in \Xs \mid \gamma_i = 0\}$.
        \label{lemma:IBC_enumeration_chemins_diriges:2}
    \end{enumerate}
\end{restatable}

\paragraph{Fork collider-free backdoor paths in \texorpdfstring{$\NC$}{NC}}\label{subsubsection:CfbF}

We first introduce an accessibility concept essential to the enumeration of fork paths. 

\begin{definition}[$\mathcal{NC}$-accessibility]
    \label{def:F-E-acc}
    Let $\Gs = (\mathcal{V}^s, \mathcal{E}^s)$ be an SCG, $P(y_t \mid \text{do}(x^1_{t-\gamma_1}), \dots, \text{do}(x^n_{t-\gamma_n}))$ be the considered effect and $V_{t_v}\in \mathcal{V}^f$.  We say that \emph{$F_{t_1}\in \mathcal{V}^f \setminus \{V_{t_v}\}$ is $V_{t_v}$-$\NC$-accessible} if there exists a candidate FTCG which contains a directed path from $F_{t_1}$ to $V_{t_v}$ which remains in $\NC$ except perhaps for $V_{t_v}$.
    We denote 
    $$t^{\NC}_{V_{t_v}}(F) \coloneqq \max \{ t_1 \mid F_{t_1} \text{is $V_{t_v}$-$\NC$-accessible} \},$$
    with the convention $\max\{\emptyset\} = -\infty$.
\end{definition}

\begin{example}
In both scenarios, \(L_t\) is \(\mathrm{Of}_t\)-\(\NC\)-accessible, since there exists a candidate FTCG containing the path \(L_t \to \mathrm{Of}_t\). Although there is also a candidate FTCG containing the path \(B_t \to L_t \to \mathrm{Of}_t\), \(B_t\) is \(\mathrm{Of}_t\)-\(\NC\)-accessible only in the second scenario, because in the first scenario \(B_t\) is itself an intervention. Consequently, in both scenarios, we have \(t^{\NC}_{\mathrm{Of}_t}(\mathrm{Living\ Room}) = t\). However,  \(t^{\NC}_{\mathrm{Of}_t}(\mathrm{Bathroom}) = -\infty\) in Scenario 1 and \(t^{\NC}_{\mathrm{Of}_t}(\mathrm{Bathroom}) = t\) in Scenario 2.
\end{example}

This leads to characterize  efficiently the existence of a collider-free backdoor path with a fork that remains in $\NC$, as proposed in  Lemma \ref{lemma:equiv_existence_chemin_fork_CD_without_consistency_through_time}.

\begin{restatable}{lemma}{mylemmaequivexistencecheminforkCDwithoutconsistencythroughtime}{(Characterization of collider-free backdoor paths with fork)}
\label{lemma:equiv_existence_chemin_fork_CD_without_consistency_through_time}
   Let $\Gs = (\mathcal{V}^s, \mathcal{E}^s)$ be an SCG and $P(y_t \mid \text{do}(x^1_{t-\gamma_1}), \dots, \Do(x^n_{t-\gamma_n}))$ be the considered effect such that for all $\Gf$ belonging to $\C(\Gs)$, $\Gf$ does not contain {a directed path from $Y_t$ to an intervention} $X^i_{t- \gamma_i} \leftsquigarrow Y_t$ which remains in $\NC \cup \{ X^i_{t-\gamma_i} \}$. The following statements are equivalent:
    \begin{enumerate}
        \item  There exists an intervention $X^i_{t-\gamma_i}$, $F_{t_f} \in \mathcal{V}^f$ and a candidate FTCG $\Gf \in \C(\Gs)$  which contains the path $X^i_{t- \gamma_i} \leftsquigarrow F_{t_f} \rightsquigarrow Y_t$ which remains in $\NC \cup \{ X^i_{t-\gamma_i} \}$.
        \label{lemma:equiv_existence_chemin_fork_CD_without_consistency_through_time:1}

        \item There exists an intervention $X^i_{t- \gamma_i}$ and $F_{t_f} \in \mathcal{V}^f$ such that $F_{t_f}$ is $X^i_{t- \gamma_i}$-$\NC$-accessible and $Y_t$-$\NC$-accessible.
        \label{lemma:equiv_existence_chemin_fork_CD_without_consistency_through_time:2}
    \end{enumerate}
\end{restatable}

The intuition behind Lemma \ref{lemma:equiv_existence_chemin_fork_CD_without_consistency_through_time} relies on a divide-and-conquer strategy to avoid searching for collider-free backdoor paths with a fork directly. In an FTCG, any such path decomposes into two directed subpaths. The lemma shows that it suffices to exhibit the first subpath in one candidate FTCG and the second subpath in another. This ensures that some candidate FTCG realizes the entire fork path, without resorting to an explicit reconstruction argument. Since testing for a single directed path in a candidate FTCG can be done efficiently, this reduction renders the overall existence check more tractable. The formal proof appears in the Supplementary Material.

Condition 2 in Lemma~\ref{lemma:equiv_existence_chemin_fork_CD_without_consistency_through_time} is not efficiently tractable, since it requires checking each \(F_{t_f} \in \mathcal{V}^f\), and the set \(\mathcal{V}^f\) is infinite. Fortunately, the set \(\{ t_1 \mid F_{t_1} \text{ is } V_{t_v}\text{-}\NC\text{-accessible} \}\) is bounded by \(t^{\NC}_{V_{t_v}}(F)\) and \(t_{\NC}(F)\) (see Lemma~\ref{lemma:CharactVNCAcc} in the appendix), and testing only the single time point \(F_{t_{\NC}(F)}\) for each time series \(F\) is sufficient (see Corollary~\ref{cor:1} in the appendix).
Consequently, Condition 2 in Lemma \ref{lemma:equiv_existence_chemin_fork_CD_without_consistency_through_time} reduces to the following:
\begin{itemize}
    \item There exists an intervention $X^i_{t- \gamma_i}$ and a time series $F$ such that $t_{\NC}(F) \leq t^{\NC}_{X^i_{t - \gamma_i}}(F)$ and $t_{\NC}(F) \leq t^{\NC}_{Y_t}(F)$,
\end{itemize}
improving further the tractability of the check of existence of a collider-free backdoor path with a fork.



\hiddensubsubsection{An Efficient Algorithm \label{sssct:algo_IBC}}

The above results show that identifiability by common adjustment in $\Gs$ is equivalent to the following two conditions:
\begin{enumerate}
    \item There does not exist an intervention $X^i_{t-\gamma_i}$ such that $ \gamma_i = 0$ and $X^i \in \Desc \left(Y, \Gs_{\mid \mathcal{S}} \right)$, where $\mathcal{S} \coloneqq \{S \in \mathcal{V}^s \mid t_{\NC}(S) \leq t \}\cup \{X^i \mid \gamma_i = 0\}$ (see Lemma~\ref{lemma:IBC_enumeration_chemins_diriges}),
    
    \item And, there does not exist an intervention $X^i_{t-\gamma_i}$ and a time series $F \in \mathcal{V}^s$ such that $t_{\NC}(F) \leq t^{\NC}_{X^i_{t - \gamma_i}}(F)$ and $t_{\NC}(F) \leq t^{\NC}_{Y_t}(F)$ (see discussion below Lemma~\ref{lemma:equiv_existence_chemin_fork_CD_without_consistency_through_time})
\end{enumerate}

Condition~\ref{lemma:IBC_enumeration_chemins_diriges:2} from Lemma~\ref{lemma:IBC_enumeration_chemins_diriges} can be verified efficiently. Indeed, it suffices to compute the set of all descendants of \(Y\) in  \(\Gs_{\mid \mathcal{S}}\) using a single breadth- or depth-first search in time \(\mathcal{O}(|\mathcal{V}^s| + |\mathcal{E}^s|)\) \citep[Chapter~22]{Cormen}, and then check if there exists an intervention \(X^i_{t-\gamma_i}\) with $\gamma_i = 0$ such that \(X^i\in\mathrm{Desc}_{\Gs_{\mid \mathcal{S}}}(Y)\).

Having ruled out directed‐paths, we now focus on the fork-path condition. To this end, we present Algorithm~\ref{algo:calcul_V_E_acc}, which computes  $\{t^{\NC}_{V_{t_v}}(F) \mid F \in \mathcal{V}^S \}$ in pseudo-linear time (see Lemma~\ref{lemma:calcul_V_E_acc} in the appendix).

As a result, the second statement of Lemma~\ref{lemma:equiv_existence_chemin_fork_CD_without_consistency_through_time} can be checked by executing Algorithm~\ref{algo:calcul_V_E_acc} twice, knowing that Algorithm \ref{algo:calcul_t_NC} has already been run. Consequently, the overall complexity is $\mathcal{O}\left(\left| \Xf \right| \log \left| \Xf \right| + (\left| \mathcal{E}^s \right| + \left| \mathcal{V}^s \right|) \log \left| \mathcal{V}^s \right|\right)$.

\begin{algorithm}[t]
\setlength{\rightskip}{-.5cm} 
\caption{Computation of $(t^{\NC}_{V_{t_v}}(S))_{S \in \mathcal{V}^s}$}\label{algo:calcul_V_E_acc}
\SetKwInput{KwData}{Input}
\SetKwInput{KwResult}{Output}
\KwData{ $\Gs = (\mathcal{V}^s, \mathcal{E}^s)$ an SCG, $\Xf$ and $V_{t_v} \in \mathcal{V}^f$}
\KwResult{$(t^{\NC}_{V_{t_v}}(S))_{S \in \mathcal{V}^s}$}
$Q \gets \text{PriorityQueue}(V_{t_v})$ \;
$t^{\NC}_{V_{t_v}}(S) \gets -\infty \quad \forall S \in \mathcal{V}^s$\;
$S.seen \gets False \quad \forall S \in \mathcal{V}^s$\;
\While{$Q \neq \emptyset$}{
    $S_{t_s} \gets Q.\text{pop\_element\_with\_max\_time\_index}()$\;
    \ForEach{unseen $P \in \text{Pa}(S, \Gs)$}{
        $t^{\NC}_{V_{t_v}}(P) \gets \max \{t_1 \mid t_1 \leq t_s \text{ and } P_{t_1} \in \NC \setminus \{V_{t_v}\}\}$\;
        \lIf{$t^{\NC}_{V_{t_v}}(P) \neq -\infty$}{
            $Q.insert(P_{t^{\NC}_{V_{t_v}}(P)})$
        }
        $P.seen \gets true$ \;
    }
}
\end{algorithm}

By combining the previous results as in Algorithm \ref{algo:calcul_IBC}, one can directly assess whether the effect is identifiable or not:

\begin{restatable}{theorem}{mythforalgoIBC}{}
    \label{th:th_for_algo_2}
    Let $\Gs = (\mathcal{V}^s, \mathcal{E}^s)$ be an SCG and $P(y_t \mid \Do(x^1_{t-\gamma_1}), \dots, \Do(x^n_{t-\gamma_n}))$ be the considered effect. Then the two statements are equivalent:
    \begin{itemize}
        \item The effect is identifiable by common adjustment in $\Gs$.
        \item Algorithm \ref{algo:calcul_IBC} outputs True.
    \end{itemize}
    Moreover Algorithm~\ref{algo:calcul_IBC} has a polynomial complexity of $\mathcal{O}\left(\left| \Xf \right|\left( \log \left| \Xf \right|+ \left( \left| \mathcal{E}^s \right| + \left| \mathcal{V}^s \right| \right)\log \left| \mathcal{V}^s \right|\right)\right)$. 
\end{restatable}

\begin{algorithm}[t]
\caption{Identifiability by common adjustment.}
\label{algo:calcul_IBC}

\SetKwInput{KwData}{Input}
\SetKwInput{KwResult}{Output}
\SetKwFunction{FMain}{AccessibleSet}
\KwData{ $\Gs = (\mathcal{V}^s, \mathcal{E}^s)$ an SCG and $\Xf$.}
\KwResult{A boolean indicating whether the effect is identifiable by common adjustment or not.}
$(t_{\NC}(S))_{S \in \mathcal{V}^s} \gets $ Algorithm \ref{algo:calcul_t_NC} \;
\tcp{Enumeration of directed paths.} 
$\mathcal{S} \gets \{S \in \mathcal{V}^s \mid t_{\NC}(S) \leq 0 \}\cup \{X^i \mid t-\gamma_i = 0\}$ \;
\If{$\exists i \in \{1,\ldots, n\}$ s.t. $ X^i \in \Desc \left(Y, \Gs_{\mid \mathcal{S}} \right)$ and $\gamma_i = 0$}{
    \Return{False}
}

\tcp{Enumeration of fork paths.} 
\ForEach{$V_{t_v} \in \{ Y_t, X^1_{t-\gamma_1}, \cdots,X^n_{t-\gamma_n} \}$}{
    $(t^{\NC}_{V_{t_v}}(S))_{S \in \mathcal{V}^s} \gets$ Algorithm \ref{algo:calcul_V_E_acc} \;
}

\ForEach{$F \in \mathcal{V}^s, X^i_{t- \gamma_i} \in (X^j_{t- \gamma_j})_j$}{
    \If{$t_{\NC}(F) \leq t^{\NC}_{X^i_{t - \gamma_i}}(F)$ and $t_{\NC}(F) \leq t^{\NC}_{Y_t}(F)$}{
        \Return{False}
    }
}

\Return{True}
\end{algorithm}

The complexity of Algorithm~\ref{algo:calcul_IBC} can be further reduced to pseudo-linear time, as detailed in Section~\ref{sct:speed_up_algo_3} of the Supplementary Material.
%
There is little interest in replacing the efficient implementation of Algorithm \ref{algo:calcul_IBC} with a formula.\footnote{We refer to a formula as a condition involving a combination of descendant, ancestor, and cycle sets, \dots} Indeed, we can not expect having a complexity better than $\mathcal{O}\left( \left| \Xf \right| + \left| \mathcal{E}^s \right| + \left| \mathcal{V}^s \right| \right)$ because in the worst case, it is necessary to traverse $\Gs$ and consider all interventions.

\hiddensection{With Consistency Through Time}\label{sec:Consistency_Time}
In practice, it is usually impossible to work with general FTCGs in which causal relations may change from one time instant to another, and people have resorted to the \textit{consistency through time} assumption (also referred to as Causal Stationarity in \citet{Runge_2018}), to obtain a simpler class of FTCGs. 

\begin{assumption}[\textit{Consistency through time}]
	\label{ass:Consistency_Time}
An FTCG $\G^f$ is said to be \emph{consistent through time} if all the causal relationships remain constant in direction through time. 
\end{assumption}

Under this assumption, the number of candidate FTCGs for a fixed SCG $\mathcal{G}^s$ is smaller, meaning that conditions to be identifiable are weaker and thus that more effects should be identifiable. We detail in  Section  \ref{sec:ctt:id} necessary and sufficient conditions to be identifiable.
All the proofs are deferred to Section \ref{sec:proof:5} in the Supplementary Material.

\hiddensubsection{Identifiability}
\label{sec:ctt:id}

Theorem \ref{th:equiv_IBC_multivarie} remains valid under Assumption \ref{ass:Consistency_Time}. Lemma \ref{lemma:IBC_enumeration_chemins_diriges} also holds because Assumption \ref{ass:Consistency_Time} only affects paths that traverse different time indices. The enumeration of collider-free backdoor paths containing a fork that remains within $\NC$, except perhaps at their first vertices, is however more complex, as detailed below. 

\begin{restatable}{lemma}{mylemmaequivexistencecheminforkNC}{}
\label{lemma:equiv_existence_chemin_fork_NC}
    Let $\Gs = (\mathcal{V}^s, \mathcal{E}^s)$ be an SCG and $P(y_t \mid \text{do}(x^1_{t-\gamma_1}), \dots, \Do(x^n_{t-\gamma_n}))$ be the considered effect such that for all $\Gf$ belonging to $\C(\Gs)$, $\Gf$ does not contain {a directed path from $Y_t$ to an intervention} $X^i_{t- \gamma_i} \leftsquigarrow Y_t$ which remains in $\NC \cup \{ X^i_{t- \gamma_i}\}$. The following statements are equivalent:
    \begin{enumerate}
        \item There exist an intervention $X^i_{t- \gamma_i}$,  $F_{t'} \in \mathcal{V}^f$ and an FTCG $\Gf \in \C(\Gs)$  containing the path $X^i_{t- \gamma_i} \leftsquigarrow F_{t'} \rightsquigarrow Y_t$ which remains in $\NC \cup \{ X^i_{t- \gamma_i}\}$.
        \label{lemma:equiv_existence_chemin_fork_NC:1}

        \item At least one of the following conditions is satisfied:
        \begin{enumerate}
            \item There exist {an intervention} $X^i_{t- \gamma_i}$ and  $F \in \mathcal{V}^s$ such that 
            $F_{t_{\NC}(F)}$ is well defined, $X^i_{t- \gamma_i}$-$\NC$-accessible and $Y_t$-$\NC$-accessible, and
            $\left\{
            \begin{aligned}
                &F \neq Y \text{, or} \\
                &t - \gamma_i \neq t_{\NC(F)}.\\
            \end{aligned}
            \right.$
            \label{lemma:equiv_existence_chemin_fork_NC:2a}
            
            \item There exists {an intervention} $X^i_{t- \gamma_i}$ such that $t - \gamma_i = t_{\NC(Y)}$ and at least one of the following properties is satisfied:
                    \begin{enumerate}
                        \item $Y_{t_{\NC}(Y)}$ is $X^i_{t- \gamma_i}$ - $\NC$-accessible without using $X^i_{t- \gamma_i} \leftarrow Y_{t- \gamma_i}$ and $Y_t$ - $\NC$-accessible.
                        \label{lemma:equiv_existence_chemin_fork_NC:2b:i}
                        
                        \item $Y_{t_{\NC}(Y)}$ is $X^i_{t- \gamma_i}$ - $\NC$-accessible and $Y_t$ - $\NC$-accessible without using $X^i_t \rightarrow Y_t$.
                        \label{lemma:equiv_existence_chemin_fork_NC:2b:ii}
                    \end{enumerate}
            \label{lemma:equiv_existence_chemin_fork_NC:2b}
        \end{enumerate}
        \label{lemma:equiv_existence_chemin_fork_NC:2}
    \end{enumerate}
\end{restatable}

Lemma \ref{lemma:equiv_existence_chemin_fork_NC} characterizes the existence of a collider-free backdoor path containing a fork. While the conditions outlined are more complex than those in Corollary~\ref{cor:1}, they play the same role, and still require only a small number of calls to $\NC$-accessibility. Consequently, one can replace the conditions in the final loop of Algorithm~\ref{algo:calcul_IBC} with conditions 2.(a) and 2.(b) of Lemma~\ref{lemma:equiv_existence_chemin_fork_NC} to derive an algorithm for identifiability by common adjustment in $\Gs$ under \textit{consistency through time}, as stated in the following theorem which is the counterpart of Theorem~\ref{th:th_for_algo_2}.

\begin{restatable}{theorem}{mythforalgoIBCwith_consistency}{}
    \label{th:th_for_algo_3}
   Let $\Gs = (\mathcal{V}^s, \mathcal{E}^s)$ be an SCG that satisfies Assumption \ref{ass:Consistency_Time} and $P(y_t \mid \Do(x^1_{t-\gamma_1}), \dots, \Do(x^n_{t-\gamma_n}))$ be the considered effect. Then the two statements are equivalent:
    \begin{itemize}
        \item The effect is identifiable by common adjustment in $\Gs$.

        \item An adaptation of Algorithm \ref{algo:calcul_IBC} outputs True.
    \end{itemize}

    In that case, a common adjustment set is given by $\C \coloneqq \left(\mathcal{V}^f \setminus \NC \right) \setminus X^f$.

\end{restatable}

In its simpler form, the adaptation of Algorithm~\ref{algo:calcul_IBC} still has a polynomial complexity of $\mathcal{O}\left(\left| \Xf \right| \cdot (\left| \mathcal{E}^s \right| + \left| \mathcal{V}^s \right| \log \left| \mathcal{V}^s \right|)\right)$. A pseudo-linear algorithm is discussed in Appendix~\ref{sct:IBC_pseudo_lineaire_consistent}.
    
The main difference between the two algorithms (with and without consistency through time) lies in how they test for collider-free backdoor paths with forks. Without consistency through time, this check is based on Lemma~\ref{lemma:equiv_existence_chemin_fork_CD_without_consistency_through_time}, while with consistency through time, this check relies on Lemma~\ref{lemma:equiv_existence_chemin_fork_NC}.
Note that assuming consistency through time reduces the number of candidate FTCGs: any candidate FTCG under consistency through time is also candidate without this assumption. As a result, the algorithm in Theorem~\ref{th:th_for_algo_2} is sound (but not complete) under consistency through time, while the algorithm from Theorem~\ref{th:th_for_algo_3} is complete (but not sound) without this assumption. 
\hiddensubsection{Experimental Illustration}

Although this article is primarily theoretical, we have conducted experiments to demonstrate the practical relevance of the results in terms of computation time and estimation. The Python implementation is available at \href{https://gricad-gitlab.univ-grenoble-alpes.fr/yvernesc/multivariateicainscg}{this repository}.\footnote{https://gricad-gitlab.univ-grenoble-alpes.fr/yvernesc/multivariateicainscg}

\begin{figure}[t]
    \vspace{14.5pt}
    \centering
    \includegraphics[width=0.4\textwidth, trim = 0 0 0 1.2cm, clip=TRUE]{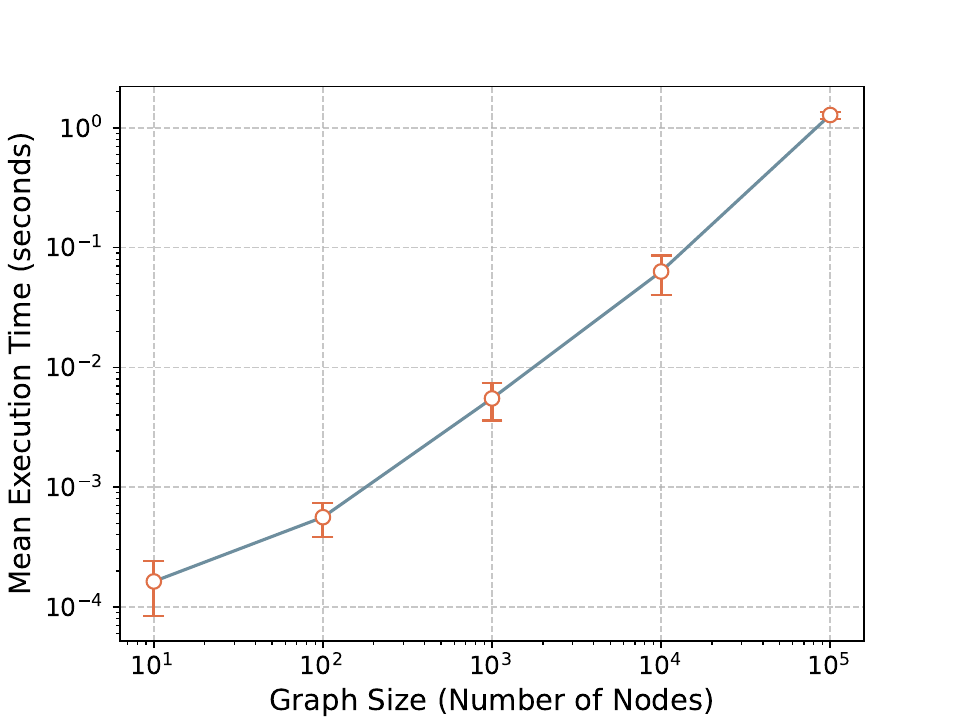}
    \caption{Average execution time of the implementation (in seconds) as a function of the number of vertices in the graph, with error bars representing standard deviation over 5 runs.}
    \label{fig:plot:execution_time}
\end{figure}
\paragraph{Execution time}
 Algorithm \ref{algo:calcul_IBC} has been implemented in Python with some speed ups discussed in Appendix \ref{sct:speed_up_algo_3}. We measure its execution speed, on a standard laptop, as a function of the graph size. For each graph size, 20 random SCGs are generated, and 5 interventions are selected at random. The average execution time of the algorithm is then measured over 5 runs and presented in Figure \ref{fig:plot:execution_time}. As one can note, even for very large graphs (with up to 100,000 vertices), the execution time remains reasonable, around 1 second, showing that the theoretical complexity of the algorithm translates into an acceptable computation time.

\begin{figure}
    \centering
    \includegraphics[width=0.4\textwidth]{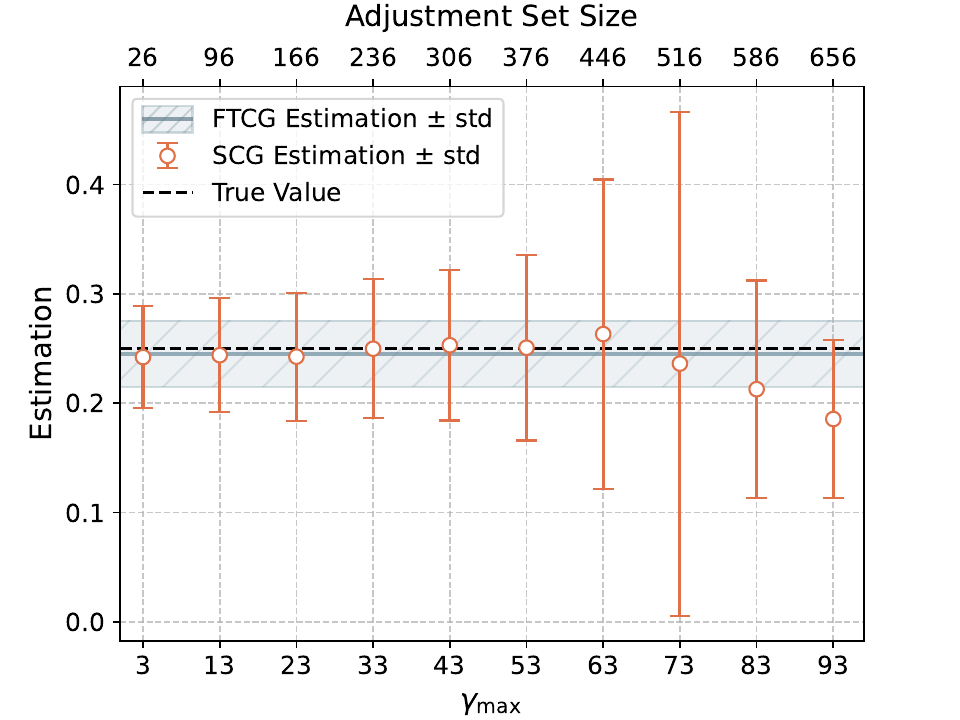}
    \caption{Estimated total causal effect \( P(y_t \mid \Do (x_t)) \) using the adjustment set from the SCG, as a function of \( \gamma_{\max} \) and as a function of the size of the adjustment set. Error bars represent the standard deviation over 100 repetitions. As a baseline, we provide the true total effect (back dashed) and the one estimated by the backdoor formulae when knowing the true FTCG, with standard deviation (gray hatched).}
    \label{fig:plot:estimations}
\end{figure}

\paragraph{Estimation}
We further considered a fixed FTCG under a linear Structural Causal Model (SCM) with additive standard Gaussian noises with a lag of 1. We designed the SCM so that the total effect is 0.25. For each choice of \(\gamma_{\max}\), which defines the farthest time horizon up to which past information is considered, we estimated the total causal effect \( P(y_t \mid \Do(x_t)) \) by using the adjustment set given by Theorem \ref{th:equiv_IBC_multivarie} up to the time index \( t - \gamma_{\max} \) over 500 data points (non overlapping windows). This estimation procedure was repeated 100 times for each \(\gamma_{\max}\). The estimated total effect and its standard deviation across these repetitions are given in Figure \ref{fig:plot:estimations}.

For comparison, we estimated the total effect using a backdoor set from the true FTCG, yielding an estimate of \(0.245\) with a standard deviation of \(0.036\) over 100 runs. The bias and variance of the SCG-based estimator remain comparable to those of the FTCG-based estimator for \(\gamma \le 50\). However, as the adjustment set expands already to 376 variables at \(\gamma = 53\), variance increases beyond large \(\gamma\), and for even larger \(\gamma\) a non-negligible bias emerges. This suggests that the adjustment set proposed in this work is particularly useful when one can assume a maximal latency of reasonable size.

\hiddensection{Conclusion}
\label{sec:conclusion}

This work has established complete conditions for identifying the effect of multiple interventions on single effects by common adjustment in summary causal graphs of time series, a criterion we have shown is both sound and complete in this setting. Specifically, Theorem \ref{th:equiv_IBC_multivarie} shows that, both with and without \textit{consistency through time}, the problem reduces to testing the existence of collider-free backdoor paths from interventions to effects that remain within a specified set. This complete characterization allowed us to derive {efficient} algorithms to determine whether an effect is identifiable by common adjustment, again both with and without \textit{consistency through time}. All the provided proofs are constructive, meaning that whenever an effect is not identifiable by common backdoor, it is possible to explicitly exhibit a collider-free backdoor path that remains within $\NC$. {Future work will focus on multiple effects and completeness results for global identifiability, \textit{i.e.}, not restricted to adjustment.} 

\begin{acknowledgements} 
    
This work was supported by fundings from the French government, managed by the National Research Agency (ANR), under the France 2030 program, reference ANR-23-IACL-0006 and  ANR-23-PEIA-0006. We thank Charles Assaad for fruitful discussions. 
\end{acknowledgements}

\newpage
\bibliography{References}

\appendix
\newpage

\onecolumn

\title{Complete Characterization for Adjustment in Summary Causal Graphs of Time Series\\(Supplementary Material)}
\maketitle

\tableofcontents

\newpage

\section{Notations}

\paragraph{Notation $\G \models$.}
When the underlying graph is not obvious from the context, we write $\G \models \phi$ to indicate that the graphical property $\phi$ holds in the graph $\G$. For example, $\G \models X \rightarrow Y$ means that $X$ is a parent of $Y$ in $\G$.

\section{Proofs of Section \ref{sec:notions}}

\mypropadjustmentSCG*

\begin{proof}
    Let us prove the two implications.
    \begin{itemize}
        \item (i) $\Rightarrow$ (ii): If $\mathbf{Z}$ satisfies the common adjustment criterion relative to $\mathbf{X}$ and $\mathbf{Y}$. Then for all $\Gf \in \C(\Gs)$, $\mathbf{Z}$ satisfies the adjustment criterion relative to $\mathbf{X}$ and $\mathbf{Y}$ in $\Gf$. Therefore, by Theorem 56 from \citep{Perkovic_2016}, Equation  \eqref{eq:adjustment} holds. 

        \item (ii) $\Rightarrow$ (i): Let us show the contrapositive.  If $\mathbf{Z}$ does not satisfy the common adjustment criterion relative to $\mathbf{X}$ and $\mathbf{Y}$. Then there exists $\Gf \in \C(\Gs)$ such that $\mathbf{Z}$ does not satisfies the adjustment criterion relative to $\mathbf{X}$ and $\mathbf{Y}$ in $\Gf$. By Theorem 57 from \citep{Perkovic_2016}, there exists $P$ compatible with $\Gf$ such that Equation \eqref{eq:adjustment} does not hold. By definition $P$ is compatible with $\Gs$ and Equation  \eqref{eq:adjustment} does not hold.
    \end{itemize}
\end{proof}

\section{How to Build an FTCG Belonging to \texorpdfstring{$\C(\Gs)$}{C(Gs)} ?}

Many proofs in this paper use constructive arguments to demonstrate the existence of an FTCG belonging to $\C(\Gs)$ that contains a given structure $\pi^f$\footnote{In most of the proofs $\pi^f$ is a path.}. To facilitate the understanding of these arguments, Lemma \ref{lemma:helper:construct_FTCG_with_path} shows how these FTCGs are constructed by adding missing arrows.

\subsection{Without Assumption \ref{ass:Consistency_Time}}

\begin{lemma}
    \label{lemma:helper:construct_FTCG_with_path}
    Let $\Gs$ be an SCG and $\pi^f$ be a graph over $\mathcal{V}^f$. If $\pi^f$ is a DAG, all its arrows respect time orientation and its reduction is a subgraph of $\Gs$, then there exists an FTCG $\Gf$ in $\C(\Gs)$ which contains $\pi^f$. 
\end{lemma}

\begin{proof}
    $\Gs$ is a SCG, by definition it is the reduction of an FTCG $\Gf_{\star} =(\mathcal{V}^f_{\star},\mathcal{E}^f_{\star})$. Let $t_{\min}$ be the minimum time index seen by $\pi^f$. Let us consider the graph $\Gf = (\mathcal{V}^f \coloneqq \mathcal{V}^f_{\star}, \mathcal{E}^f)$ whose edges are constructed as follows:
    \begin{enumerate}
        \item All edges from $\pi^f$ are set in $\mathcal{E}^f$.
        \label{lemma:helper:construct_FTCG_with_path:step1}
        
        \item For all the edges $A \rightarrow B$ in $\Gs$ that are not reductions of arrows in $\pi^f$, add the edge $A_{t_{\min} -1} \rightarrow B_{t}$ into $\mathcal{E}^f$. 
        \label{lemma:helper:construct_FTCG_with_path:step2}

    \end{enumerate}
    
    We have the following properties:
    \begin{itemize}

        \item $\Gf$ contains $\pi^f$ because all arrows of $\pi^f$ are set in $\mathcal{E}^f$.
                
        \item $\Gf$ is a DAG. Indeed, there is no cycle due to instantaneous arrows because all instantaneous arrows come  from $\pi^f$ and $\pi^f$ is a DAG. Moreover, there is no cycle due to delayed arrows because all of them follow the flow of time.
                
        \item $\Gf$ belongs to $\C(\Gs)$: Let us consider $\G^r = (\mathcal{V}^r, \mathcal{E}^r)$ the reduction of $\Gf$. By definition, $\mathcal{V}^r = \mathcal{V} = \mathcal{V}^s$ . Les us prove that $\mathcal{E}^r = \mathcal{E}^s$ by showing the two inclusions:
        \begin{itemize}
            \item $\mathcal{E}^r \subseteq \mathcal{E}^s$: Let us consider an arrow $a^r$ of $G^r$ and $a^f$ a corresponding arrow in $\Gf$. We distinguish two cases:
            \begin{itemize}
                \item If $a^f$ is in $\pi^f$ then $a^r$ is in $\pi^r$ and $\pi^r$ is a subgraph of $\Gs$. Therefore  $a^r \in \mathcal{E}^s$.

                \item Otherwise, $a^r$ have been added to $\mathcal{E}^f$ during step \ref{lemma:helper:construct_FTCG_with_path:step2}. Therefore $a^r \in \mathcal{E}^s$.
            \end{itemize}
            In both cases $a^r \in \mathcal{E}^s$, therefore $\mathcal{E}^r \subseteq \mathcal{E}^s$.
            
            \item $\mathcal{E}^s \subseteq \mathcal{E}^r$:  Let us consider an arrow $A \rightarrow B$ of $G^s$. we distinguish two cases:
            \begin{itemize}
                \item If $A \rightarrow B$ is the reduction of an arrow $a^f$ of $\pi^f$, then $a^f \in \mathcal{E}^f$. Therefore $A \rightarrow B$ is an arrow of $\G^r$

                \item Otherwise, $A_{t_{\min} -1} \rightarrow B_{t}$ is in $\mathcal{E}^f$. Therefore $A \rightarrow B$ is an arrow of $\G^r$
            \end{itemize} 
            In both cases $A \rightarrow B$ is an arrow of $\G^r$, therefore $\mathcal{E}^s \subseteq \mathcal{E}^r$.
        \end{itemize}  
        Therefore, $\G^r = \G^s$ and thus $\Gf$ belongs to $\C(\Gs)$.
    \end{itemize}
    Therefore, $\Gf$ is an FTCG belonging to $\C(\Gs)$ which contains $\pi^f$.
\end{proof}

\subsection{With Assumption \ref{ass:Consistency_Time}}

The construction given in Lemma \ref{lemma:helper:construct_FTCG_with_path} does not work to build an FTCG that satisfies Assumption \ref{ass:Consistency_Time}. Indeed, some arrows are missing. By copying the construction of lemma \ref{lemma:helper:construct_FTCG_with_path} at each time step, if $\pi^f$ verify the correct properties, it is possible to construct an FTCG belonging to $\C(\Gs)$ which contains $\pi^f$ and satisfy assumption \ref{ass:Consistency_Time}. This reasoning is encapsulated by Lemmma \ref{lemma:helper:construct_FTCG_with_path_with_cst}.

\begin{lemma}
    \label{lemma:helper:construct_FTCG_with_path_with_cst}
    Let $\Gs$ be an SCG and $\pi^f$ be a graph over $\mathcal{V}^f$. If $\pi^f$ is a DAG, all its arrows respect time orientation, its reduction is a subgraph of $\Gs$ and all its instantaneous arrows respect Assumption \ref{ass:Consistency_Time}, then there exists an FTCG $\Gf$ in $\C(\Gs)$ which contains $\pi^f$ and satisfies Assumption \ref{ass:Consistency_Time}. 
\end{lemma}

\begin{proof}
    $\Gs$ is a SCG, by definition it is the reduction of an FTCG $\Gf_{\star} =(\mathcal{V}^f_{\star},\mathcal{E}^f_{\star})$. Let $t_{\min}$ be the minimum time index seen by $\pi^f$. Let us consider the graph $\Gf = (\mathcal{V}^f \coloneqq \mathcal{V}^f_{\star}, \mathcal{E}^f)$ whose edges are constructed as follows:
    \begin{enumerate}
        \item All edges from $\pi^f$ are set in $\mathcal{E}^f$.
        \label{lemma:helper:construct_FTCG_with_path_with_cst:step1}
        
        \item For all the edges $A \rightarrow B$ in $\Gs$ that are not reductions of arrows in $\pi^f$, add the edge $A_{t_{\min} -1} \rightarrow B_{t}$ into $\mathcal{E}^f$. 
        \label{lemma:helper:construct_FTCG_with_path_with_cst:step2}

        \item Copy these arrows at each time step.
        \label{lemma:helper:construct_FTCG_with_path_with_cst:step3}
    \end{enumerate}

    We have the following properties:
    \begin{itemize}

        \item $\Gf$ contains $\pi^f$ because all arrows of $\pi^f$ are set in $\mathcal{E}^f$.
                
        \item $\Gf$ is a DAG: All instantaneous arrows come from $\pi^f$ and its copies. Since $\pi^f$ is a DAG and all its instantaneous arrows respect Assumption \ref{ass:Consistency_Time}, we know that there is no cycle made of instantaneous arrows. Moreover, delayed arrows follow the flow of time, therefore, they cannot form a cycle. Therefore, $\Gf$ does not contain any cycle.
        
        Moreover, there is no cycle due to delayed arrows because all of them follow the flow of time.
                
        \item $\Gf$ belongs to $\C(\Gs)$: By the same reasoning of Lemma \ref{lemma:helper:construct_FTCG_with_path}, at Step \ref{lemma:helper:construct_FTCG_with_path_with_cst:step2}, the reduction $\G^r_{\text{at Step \ref{lemma:helper:construct_FTCG_with_path_with_cst:step2}}}$ of $\G^f_{\text{at Step \ref{lemma:helper:construct_FTCG_with_path_with_cst:step2}}}$ is equal to $\G^s$. Adding copies of arrows does not change the reduction. Therefore, the reduction of $\Gf$ is equal to $\Gs$ i.e. $\Gf$ belongs to $\C(\Gs)$.

        \item $\Gf$ satisfies Assumption \ref{ass:Consistency_Time} because copying all arrows ensure that all causal relationships remain constant throughout time
    \end{itemize}
    Therefore, $\Gf$ is an FTCG belonging to $\C(\Gs)$ which contains $\pi^f$ and satisfies Assumption \ref{ass:Consistency_Time}.
\end{proof}

\newpage
\section{Proofs of Section \ref{section:ICA}}\label{sec:proof:4}

\subsection{Proofs of Section \ref{subsection:equiv_based_CD}}

\subsubsection{Proof of Lemma \ref{lemma:def_equiv_CD}}

\mylemmadefequivCD*

\begin{proof}
    \textbf{We start by proving the characterizations of $\NC$.} \newline
    \noindent By definition, we have: \newline
    \[
    \NC \coloneqq \bigcup_{\Gf \in \mathcal{C}(\Gs)} \Forb\left(\Xf, Y_t, \Gf\right)\setminus \Xf
    \]

    Le us denote  $\mathcal{S} \coloneqq  \bigcup_{Z \in \mathcal{V}^S} \{Z_{t_1}\}_{t_1 \geq t_{\NC}(Z)} \setminus \Xf$. We prove $\NC = \mathcal{S}$ by double inclusion:
    \begin{itemize}
        \item $\NC \subseteq \mathcal{S}$: Let $Z_{t_1} \in \NC$. By definition, $t_1 \geq t_{\NC}(Z)$ and $Z_{t_1} \notin \Xf$.  Therefore $\NC \subseteq \mathcal{S}$.
        
        \item $\mathcal{S} \subseteq \NC$: Let $Z_{t_1} \in \mathcal{S}$. By definition $Z_{t_1} \notin \Xf$. Necessarily $t_{\NC}(Z) < + \infty$ thus there exists an FTCG $\Gf \in \C(\Gs)$ in which $Z_{t_{\NC}} \in \Forb\left(\Xf, Y_t, \Gf\right)$. Thus, in $\Gf$, there exists a vertex $W_{t_w} \notin \Xf$ which lies on a proper causal path $\pi^f_1$ from $\Xf$ to $Y_t$ and $\pi^f_2$ a directed path from $W_{t_w}$ to $Z_{t_{\NC}(Z)}$. We can construct another FTCG belonging to $\C(\Gs)$ containing $\pi^f_1$ and $\pi^f_2$ except that the last arrow of $\pi^f_2$ points to $Z_{t_1}$ instead of $Z_{t_{\NC}(Z)}$. Thus, $Z_{t_1} \in \NC$. Therefore $\mathcal{S} \subseteq \NC $.
    \end{itemize}

    \noindent Therefore, $\NC = \bigcup_{Z \in \mathcal{V}^S} \{Z_{t_1}\}_{t_1 \geq t_{\NC}(Z)} \setminus \Xf$.

\textbf{Lemma \ref{lemma:preuve_algo_t_NC} proves that Algorithm \ref{algo:calcul_t_NC} computes computes $(t_{\NC}(S))_{S \in \mathcal{V}^s}$ in pseudo-linear complexity.} The proof requires several additional lemmas. It is discussed in the following paragraphs.
\end{proof}

To compute  $(t_{\NC}(F))_{F \in \mathcal{V}^s}$ in pseudo-linear complexity, we need to use an efficient characterization of $\CF$ given by Lemma \ref{lemma:characterisation_CF}.

\begin{lemma}[Characterization of $\CF$.]
\label{lemma:characterisation_CF}
Let $\Gs = (\mathcal{V}^s, \mathcal{E}^s)$ be an SCG and let $P(y_t \mid \text{do}(x^1_{t-\gamma_1}), \dots, \text{do}(x^n_{t-\gamma_n}))$ be the considered effect. For each $C \in \Ch(\Xs)$, we define:

\[
t_C \coloneqq \min \{ t_1 \mid C_{t_1} \notin \Xf ~\text{and}~ \exists i \exists \Gf \in \C(\Gs) \text{ s.t } \Gf \models X^i_{t-\gamma_i} \rightarrow C_{t_1} \text{ and }  C_{t_1} \in \Anc(Y_t, \Gf\setminus \Xf) \},
\]

\noindent with the convention that $\min \emptyset = +\infty$. We have:
\[
\CF = \bigcup_{\substack{C ~\in~ \Ch\left(\Xs\right)  \\ t_C~ <~ +\infty}}
        \ \ \bigcup_{\substack{\Gf \in \C\left(\Gs\right) \\ \exists i ~ \Gf \models X^i_{t-\gamma_i} \rightarrow C_{t_C}}}
        \Desc\left(C_{t_C}, \Gf \right).
\]
\end{lemma}

\begin{proof}
    We prove this lemma by double inclusion. Let \( S = \bigcup_{\substack{C ~\in~ \Ch\left(\Xs\right)  \\ t_C~ <~ +\infty}}
        \ \ \bigcup_{\substack{\Gf \in \C\left(\Gs\right) \\ \exists i ~ \Gf \models X^i_{t-\gamma_i} \rightarrow C_{t_C}}}
        \Desc\left(C_{t_C}, \Gf \right)\):
    \begin{itemize}
        \item Let $D_{t_d} \in \CF$. By definition there exists an FTCG $\Gf \in \C(\Gs)$ such that $D_{t_d} \in \Forb\left(\Xf, Y_t, \Gf \right)$. Thus, $\Gf$ contains a proper causal path $\pi^f$ from $\Xf$ to $Y_t$ and there exists $W_{t_w} \in \pi^f \setminus \Xf$ such that $D_{t_d} \in \Desc(W_{t_w}, \Gf)$. Let $C_{t_1}$ be the second vertex of $\pi^f$. Let $X ^i$ be the parent of $C$ in $\Gs$ with the smallest corresponding $t-\gamma_i$. Let $t_\star$ be the minimal time greater than $t-\gamma_i$ such that $C_{t_\star}$ is not an intervention. By changing the first two vertices of $\pi^f$, we can construct an FTCG ${\Gf}'$ that contains ${\pi^f}' = X^i_{t-\gamma_i} \rightarrow C_{t_\star} \rightsquigarrow Y_t$ such that $W_{t_w} \in {\pi^f}' \setminus \Xf$ and $D_{t_d} \in \Desc(W_{t_w}, {\Gf}')$. By construction,  $C_{t_\star} \notin \Xf$, ${\Gf}' \models X^i_{t-\gamma_i} \rightarrow C_{t_\star}$ and $ C_{t_\star} \in \Anc(Y_t, {\Gf}' \setminus \Xf)$. Thus $t_C \leq t_\star < +\infty$. $X^i$ is the parent of $C$ in $\Gs$ with the smallest corresponding $t-\gamma_i$ thus $t_C = t_\star$. Thus, $W_{t_w} \in \Desc(C_{t_c}, {\Gf}')$ and we already know that $D_{t_d} \in \Desc(W_{t_w}, {\Gf}')$. Thus $D_{t_d} \in \Desc(C_{t_C}, {\Gf}')$. Therefore $D_{t_d} \in S.$

        \item Let $D_{t_d} \in S$. By definition, $t_c < +\infty$, thus there exist an intervention $X^i_{t-\gamma_i}$ and an FTCG $\Gf_1$ such that $\Gf_1 \models X^i_{t-\gamma_i} \rightarrow C_{t_C}$ and $C_{t_C} \in \Anc(Y_t, \Gf_1 \setminus \Xf)$. Let $\pi^f_1$ be this path. Without loss of generality, we can assume that $\pi^f_1$ does not pass twice through the same time series. Moreover, there exist an intervention $X^j_{t-\gamma_j}$ and an FTCG $\Gf_2$ such that $\Gf_2 \models X^j_{t-\gamma_j} \rightarrow C_{t_C} \rightsquigarrow D_{t_d}$. Let $\pi^f_2$ be this path. We construct $\Gf_3$ with the following procedure:
        \begin{enumerate}
            \item $\Gf_3$ is the graph without edges on $\mathcal{V}^f$.
            \item Insert $\pi^f_1$ and $X^i_{t-\gamma_i} \rightarrow C_{t_C}$ into $\Gf_3$. 
            \item At this point of the construction, we know that $C_{t_C}$ is on a proper causal path from $\Xf$ to $Y_t$ in $\Gf_3$.
            \item Let $W_{t_w} \in \pi^f_2$ be such that $W$ is the last time series seen by $\pi^f_2$ that is also seen by $\pi^f_1$. Let $\langle V^1_{t_1} = C_{t_C}, \dots, V^k_{t_k} \rangle$ denote the vertex in $\pi^f_1$ until it reaches $W$. For all $m$ from $1$ to $k-1$, insert $V^m_{t_C} \rightarrow V^{m+1}_{t_C}$ into $\Gf_3$. Insert $V^{k-1}_{t_C} \rightarrow W_{t_w}$ into $\Gf_3$. Insert all the arrows from $\pi^f_2$ after $W_{t_w}$.
            \item At this point of the construction, $D_{t_d}$ is a descendant of $C_{t_C}$ in $\Gf_3$.
            \item Add all missing arrow to a greater time lag so that $\Gf_3 \in \C(\Gs)$.
        \end{enumerate}
        Therefore, $D_{t_d} \in \CF$.
    \end{itemize}

    Therefore, $\CF = S$.
    One can check that indeed $\Gf_3 \in \C(\Gs)$ and that $\Gf_3$ satisfies Assumption \ref{ass:Consistency_Time} if both $\Gf_1$ and $\Gf_2$ satisfy Assumption \ref{ass:Consistency_Time}.
\end{proof}

Lemma \ref{lemma:characterisation_CF_suite} specifies Lemma \ref{lemma:characterisation_CF}.
\begin{lemma}
\label{lemma:characterisation_CF_suite}
    Let $\Gs = (\mathcal{V}^s, \mathcal{E}^s)$ be an SCG and let $P(y_t \mid \text{do}(x^1_{t-\gamma_1}), \dots, \text{do}(x^n_{t-\gamma_n}))$ be the considered effect. Let $C \in \Ch(\Xs)$ such that $t_C < +\infty$. We have:

    \[
    \bigcup_{\substack{\Gf \in \C\left(\Gs\right) \\ \exists i ~ \Gf \models X^i_{t-\gamma_i} \rightarrow C_{t_C}}}
        \Desc\left(C_{t_C}, \Gf \right) = \begin{cases}
                        \bigcup_{\substack{\Gf \in \C\left(\Gs\right)}} \Desc\left(C_{t_C}, \Gf \right) & \text{if } d(C),\\
                        \bigcup_{\substack{\Gf \in \C\left(\Gs\right)}} \Desc\left(C_{t_C}, \Gf \setminus \{X^i_{t-\gamma_i}\} \right)  & \text{otherwise.}
                    \end{cases}
    \]

Where $d(C) \coloneqq 
                     \# \{i \mid X^i \in \text{Pa}(C, \Gs) \text{ and }t-\gamma_i < t_C \} \geq 1 \text{ or }  
                     \# \{i \mid X^i \in \text{Pa}(C, \Gs) \text{ and } t-\gamma_i =t_C \} \geq 2$, and, in the second case, $X^i$ is the only element of $\{X^i \in \text{Pa}(C, \Gs) \text{ and } t-\gamma_i \leq t_C \} = \{X^i \in \text{Pa}(C, \Gs) \text{ and } t-\gamma_i = t_C \}$.
\end{lemma}
    
\begin{proof}
$t_C < +\infty$ thus \(\# \{i \mid X^i \in \text{Pa}(C, \Gs) \text{ and } t-\gamma_i \leq t_C \} \geq 1 \). We distinguish two cases:
\begin{itemize}
    \item If $\# \{i \mid X^i \in \text{Pa}(C, \Gs) \text{ and }t-\gamma_i < t_C \} \geq 1$.  Let us show the two inclusions:
    \begin{itemize}
        \item By definition, $\bigcup_{\substack{\Gf \in \C\left(\Gs\right) \\ \exists i ~ \Gf \models X^i_{t-\gamma_i} \rightarrow C_{t_C}}}
            \Desc\left(C_{t_C}, \Gf \right) \subseteq \bigcup_{\substack{\Gf \in \C\left(\Gs\right)}} \Desc\left(C_{t_C}, \Gf \right)$.

        \item Let $X^i_{t-\gamma_i}$ be an element of $\{X^i \in \text{Pa}(C, \Gs) \mid t-\gamma_i < t_C \}$.  The arrow from $X^i_{t-\gamma_i}$ to $C_{t_C}$ cannot contradict any other arrow. Therefore, $\bigcup_{\substack{\Gf \in \C\left(\Gs\right)}} \Desc\left(C_{t_C}, \Gf \right) \subseteq \bigcup_{\substack{\Gf \in \C\left(\Gs\right) \\ \exists i ~ \Gf \models X^i_{t-\gamma_i} \rightarrow C_{t_C}}} \Desc\left(C_{t_C}, \Gf \right)$.
    \end{itemize} 
    Therefore, $\bigcup_{\substack{\Gf \in \C\left(\Gs\right) \\ \exists i ~ \Gf \models X^i_{t-\gamma_i} \rightarrow C_{t_C}}} \Desc\left(C_{t_C}, \Gf \right) = \bigcup_{\substack{\Gf \in \C\left(\Gs\right)}} \Desc\left(C_{t_C}, \Gf \right)$.

    \item Otherwise, \(\# \{i \mid X^i \in \text{Pa}(C, \Gs) \text{ and }  t-\gamma_i = t_C  \} \geq 1 \). We distinguish two cases:
    \begin{itemize}
        \item  If \(\# \{i \mid X^i \in \text{Pa}(C, \Gs) \text{ and } t-\gamma_i = t_C \} = 1 \). Let $X^i_{t-\gamma_i}$ be the only element of $\{ X^i_{t-\gamma_i} \mid X^i \in \text{Pa}(C, \Gs) \mid  t-\gamma_i = t_C \}$. In this case, $\bigcup_{\substack{\Gf \in \C\left(\Gs\right) \\ \exists i ~ \Gf \models X^i_{t-\gamma_i} \rightarrow C_{t_C}}}
        \Desc\left(C_{t_C}, \Gf \right) = \bigcup_{\substack{\Gf \in \C\left(\Gs\right) \\\Gf \models X^i_{t-\gamma_i} \rightarrow C_{t-\gamma_i}}}
        \Desc\left(C_{t-\gamma_i}, \Gf \right)$. Let us show the two inclusions:
        \begin{itemize}
            \item Let $D_{t_d} \in \bigcup_{\substack{\Gf \in \C\left(\Gs\right) \\\Gf \models X^i_{t-\gamma_i} \rightarrow C_{t-\gamma_i}}}$. There exists an FTCG $\Gf$ which contains a directed path $\pi^f$ from $C_{t-\gamma_i}$ to $D_{t_d}$ and such that $\Gf \models X^i_{t-\gamma_i} \rightarrow C_{t-\gamma_i}$. $\pi^f$ does not pass through $X^i_{t-\gamma_i}$ because, otherwise, $\Gf \models X^i_{t-\gamma_i} \leftarrow C_{t-\gamma_i}$. Thus $D_{t_d} \in \bigcup_{\substack{\Gf \in \C\left(\Gs\right)}} \Desc\left(C_{t_C}, \Gf \setminus \{X^i_{t-\gamma_i}\} \right)$.

            \item Let $D_{t_d} \in \bigcup_{\substack{\Gf \in \C\left(\Gs\right)}} \Desc\left(C_{t_C}, \Gf \setminus \{X^i_{t-\gamma_i}\} \right)$. Thus, there exists an FTCG containing a directed path $\pi^f$ from $C_{t-\gamma_i}$ to $D_{t_d}$ without passing through $X^i_{t-\gamma_i}$. Without loss of generality, we can assume that $\pi^f$ does not pass twice through the same time series. We distinguish two cases:
            \begin{itemize}
                \item If $\pi^f$ does not use an arrow from $C_{t-\gamma_i}$ to some $X^i_{t'} \in X^i$. Then, we can construct an FTCG containing $\pi^f$ and $X^i_{t-\gamma_i} \rightarrow C_{t-\gamma_i}$. Thus, $D_{t_d} \in \bigcup_{\substack{\Gf \in \C\left(\Gs\right) \\\Gf \models X^i_{t-\gamma_i} \rightarrow C_{t-\gamma_i}}}$.

                \item Otherwise, $\pi^f$ does not pass through $X^i_{t-\gamma_i}$, thus $t_d < t-\gamma_i$. We can construct an FTCG containing $X^i_{t-\gamma_i} \rightarrow C_{t-\gamma_i}$ and a path ${\pi^f}'$ from $C_{t-\gamma_i}$ to $D_{t_d}$ where the first arrow is $C_{t-\gamma_i} \rightarrow X^i_{t_d}$ and all other arrows are the arrows from $\pi^f$ but at time $t_d$. Thus, $D_{t_d} \in \bigcup_{\substack{\Gf \in \C\left(\Gs\right) \\\Gf \models X^i_{t-\gamma_i} \rightarrow C_{t-\gamma_i}}}$.
            \end{itemize}
            Therefore, $D_{t_d} \in \bigcup_{\substack{\Gf \in \C\left(\Gs\right) \\\Gf \models X^i_{t-\gamma_i} \rightarrow C_{t-\gamma_i}}}$.
        \end{itemize}
        Therefore, $\bigcup_{\substack{\Gf \in \C\left(\Gs\right) \\\Gf \models X^i_{t-\gamma_i} \rightarrow C_{t-\gamma_i}}} = \bigcup_{\substack{\Gf \in \C\left(\Gs\right)}} \Desc\left(C_{t_C}, \Gf \setminus \{X^i_{t-\gamma_i}\} \right)$.

        \item Otherwise, \(\# \{i \mid X^i \in \text{Pa}(C, \Gs) \text{ and } t-\gamma_i = t_C \} \geq 2 \). Let us show the two inclusions:
        \begin{itemize}
            \item By definition, $\bigcup_{\substack{\Gf \in \C\left(\Gs\right) \\ \exists i ~ \Gf \models X^i_{t-\gamma_i} \rightarrow C_{t_C}}}\Desc\left(C_{t_C}, \Gf \right) \subseteq \bigcup_{\substack{\Gf \in \C\left(\Gs\right)}} \Desc\left(C_{t_C}, \Gf \right)$.

            \item Let $D_{t_d} \in \bigcup_{\substack{\Gf \in \C\left(\Gs\right)}} \Desc\left(C_{t_C}, \Gf \right)$. There exist an FTCG which contains a directed path $\pi^f$ from $C_{t_C}$ to $D_{t_d}$. Without loss of generality, we can assume that $\pi^f$ does not pass twice through the same time series. Thus, $\pi^f$ contains only one arrow $a$ coming from $C_{t_C}$.  \(\# \{i \mid X^i \in \text{Pa}(C, \Gs) \text{ and } t-\gamma_i = t_C \} \geq 2 \) thus there exists an intervention $X^i_{t-\gamma_i}\in \{X^i_{t-\gamma_i} \mid X^i \in \text{Pa}(C, \Gs) \text{ and } t-\gamma_i = t_C \}$ such that $X^i_{t-\gamma_i} \rightarrow C_{t_C}$ does not contradict $a$. Thus we can construct another FTCG containing $X^i_{t-\gamma_i} \rightarrow C_{t_C}$ and $\pi^f$. Thus, $\bigcup_{\substack{\Gf \in \C\left(\Gs\right)}} \Desc\left(C_{t_C}, \Gf \right) \subseteq \bigcup_{\substack{\Gf \in \C\left(\Gs\right) \\ \exists i ~ \Gf \models X^i_{t-\gamma_i} \rightarrow C_{t_C}}} \Desc\left(C_{t_C}, \Gf \right)$
        \end{itemize}
        Therefore, $\bigcup_{\substack{\Gf \in \C\left(\Gs\right) \\ \exists i ~ \Gf \models X^i_{t-\gamma_i} \rightarrow C_{t_C}}} \Desc\left(C_{t_C}, \Gf \right) = \bigcup_{\substack{\Gf \in \C\left(\Gs\right)}} \Desc\left(C_{t_C}, \Gf \right)$.
    \end{itemize}
\end{itemize}
\end{proof}

Lemm \ref{lemma:calcul_t_C} shows how to compute $t_C$.
\begin{lemma}
\label{lemma:calcul_t_C}
Let $\Gs = (\mathcal{V}^s, \mathcal{E}^s)$ be an SCG and let $P(y_t \mid \text{do}(x^1_{t-\gamma_1}), \dots, \text{do}(x^n_{t-\gamma_n}))$ be the considered effect. Let $C \in \Ch(\Xs)$. let $t_{\max} \coloneqq\max \left\{ t_1 \mid  \exists\Gf\text{ s.t }  C_{t_1} \in \Anc(Y_{t}, \Gf \setminus \Xf)  \right\}$ and $t_{\min} \coloneqq \min \{t-\gamma_i \mid X^i \in \text{Pa}(C,\Gs)\}$. We have the following identity:
\[
\left\{ t_1 \mid C_{t_1} \notin \Xf ~\text{and}~ \exists i \exists \Gf \in \C(\Gs) \text{ s.t } \Gf \models X^i_{t-\gamma_i} \rightarrow C_{t_1} \text{ and }  C_{t_1} \in \Anc(Y_t, \Gf\setminus \Xf) \right\} = 
\left\{ t_1 \in \left[ t_{\min}, t_{\max} \right]  \mid C_{t_1} \notin \Xf  \right\}.
\]
\end{lemma}
\begin{proof}
    Let us prove the two inclusions.
    \begin{itemize}
        \item Let $t_1 \in \{ t_1 \mid C_{t_1} \notin \Xf ~\text{and}~ \exists i \exists \Gf \in \C(\Gs) \text{ s.t } \Gf \models X^i_{t-\gamma_i} \rightarrow C_{t_1} \text{ and }  C_{t_1} \in \Anc(Y_t, \Gf\setminus \Xf) \}$. By definition, there exist $\Gf$ and an intervention $X^i_{t-\gamma_i}$ such that $\Gf \models X^i_{t-\gamma_i} \rightarrow C_{t_1}$. Thus $\Gs \models X^i \rightarrow C$. Thus $t_1 \geq t-\gamma_i \geq t_{\min}$ because causality does not move backwards in time. Moreover, $ C_{t_1} \in \Anc(Y_t, \Gf\setminus \Xf)$. Thus $t_1 \leq t_{\max}$. Therefore $\{ t_1 \mid C_{t_1} \notin \Xf ~\text{and}~ \exists i \exists \Gf \in \C(\Gs) \text{ s.t } \Gf \models X^i_{t-\gamma_i} \rightarrow C_{t_1} \text{ and }  C_{t_1} \in \Anc(Y_t, \Gf\setminus \Xf) \} \subseteq \left\{ t_1 \in \left[ t_{\min},t_{\max} \right]  \mid C_{t_1} \notin \Xf  \right\}$.

        \item Let $t_1 \in  \left\{ t_1 \in \left[ t_{\min}, t_{\max} \right]  \mid C_{t_1} \notin \Xf  \right\}$. By definition $t_1 \leq t_{\max}$. Thus there exist an FTCG containing $\pi^f$, a directed path from $C_{t_{\max}}$ to $Y_t$ that does not go through $\Xf$. By definition $t_1 \geq t_{\min}$, thus there exists an intervention $X^i_{t-\gamma_i}$ such that $t-\gamma_i = t_{\min} \leq t_1$ and $X^i \in \text{Pa}(C,\Gs)$. Since $\pi^f$ does not go through $\Xf$, by changing its first arrow we can construct an FTCG $\Gf \in \C(\Gs)$ such that $\Gf \models X^i_{t-\gamma_i} \rightarrow C_{t_1} \text{ and }  C_{t_1} \in \Anc(Y_t, \Gf\setminus \Xf)$. Therefore  $ \left\{ t_1 \in \left[ t_{\min},t_{\max} \right]  \mid C_{t_1} \notin \Xf  \right\} \subseteq \{ t_1 \mid C_{t_1} \notin \Xf ~\text{and}~ \exists i \exists \Gf \in \C(\Gs) \text{ s.t } \Gf \models X^i_{t-\gamma_i} \rightarrow C_{t_1} \text{ and }  C_{t_1} \in \Anc(Y_t, \Gf\setminus \Xf) \}$.
    \end{itemize}
    Therefore, $\{ t_1 \mid C_{t_1} \notin \Xf ~\text{and}~ \exists i \exists \Gf \in \C(\Gs) \text{ s.t } \Gf \models X^i_{t-\gamma_i} \rightarrow C_{t_1} \text{ and }  C_{t_1} \in \Anc(Y_t, \Gf\setminus \Xf) \} = 
\left\{ t_1 \in \left[ t_{\min}, t_{\max} \right]  \mid C_{t_1} \notin \Xf  \right\}.$
\end{proof}

\begin{lemma}
\label{lemma:truc}
 Let $\Gs = (\mathcal{V}^s, \mathcal{E}^s)$ be an SCG and let $P(y_t \mid \text{do}(x^1_{t-\gamma_1}), \dots, \text{do}(x^n_{t-\gamma_n}))$ be the considered effect. Let $C \in \Ch(\Xs)$ and $S$ be a time series. We have the following identity:

 \begin{equation}
 \label{eq:A}
 \min \left\{ t_1 \mid S_{t_1} \in \bigcup_{\substack{\Gf \in \C\left(\Gs\right)}}\Desc(C_{t_C}, \Gf)\setminus \Xf\right\}  = \begin{cases}
     \min \{t_1 \mid t_1 \geq t_C \text{ and } S_{t_1} \notin \Xf\} &\text{if } S \in \Desc(C, \Gs),\\
     +\infty & \text{otherwise.}
 \end{cases}
 \end{equation}

\noindent When $t_C = t-\gamma_i$, we also have:

\begin{equation}
 \label{eq:B}
 \min \left\{ t_1 \mid S_{t_1} \in \bigcup_{\substack{\Gf \in \C\left(\Gs\right)}}\Desc(C_{t_C}, \Gf \setminus \{ X ^i_{t-\gamma_i}\}) \setminus \Xf\right\} = \begin{cases}
     \min \{t_1 \mid t_1 \geq t_C \text{ and } S_{t_1} \notin \Xf\} &\text{if } S \in \Desc(C, \Gs\setminus \{ X^i\}),\\
     \min \{t_1 \mid t_1 \geq t_C+1 \text{ and } S_{t_1} \notin \Xf\} & \text{else if } S \in \Desc(C, \Gs ),\\
     +\infty & \text{otherwise.}
 \end{cases}
 \end{equation}
\end{lemma}

\begin{proof}

Let us prove Equation \eqref{eq:A}. We distinguish two cases:
\begin{itemize}
    \item If $S \in \Desc(C, \Gs)$, then there exists a directed path $\pi^s$ from $C$ to $S$ in $\Gs$. For all $t_1 \geq t_C$, there exists an FTCG which contains a directed path from $C_{t_C}$ to $S_{t_1}$. Therefore $\min \left\{ t_1 \mid S_{t_1} \in \bigcup_{\substack{\Gf \in \C\left(\Gs\right)}}\Desc(C_{t_C}, \Gf)\setminus \Xf\right\}  =  \min \{t_1 \mid t_1 \geq t_C \text{ and } S_{t_1} \notin \Xf\}$.

    \item Otherwise, necessarily, for all $t_1$, $S_{t_1}$ cannot be a descendant of $C_{t_C}$ in any FTCG because otherwise S would be a descendant of $C$ in $\Gs$. Thus $\min \left\{ t_1 \mid S_{t_1} \in \bigcup_{\substack{\Gf \in \C\left(\Gs\right)}}\Desc(C_{t_C}, \Gf)\setminus \Xf\right\}  =  +\infty$.
\end{itemize}

Therefore, Equation \eqref{eq:A} holds.

Let us prove Equation \eqref{eq:B}. We distinguish three cases:

\begin{itemize}
    \item If $S \in \Desc(C, \Gs \setminus \{ X^i \})$, then there exists a directed path $\pi^s$ from $C$ to $S$ in $\Gs \setminus \{ X^i \}$. For all $t_1 \geq t_C$, there exists an FTCG $\Gf$  such that $\Gf \setminus \{ X^i_{t-\gamma_i} \}$ contains a directed path from $C_{t_C}$ to $S_{t_1}$. Therefore $\min \left\{ t_1 \mid S_{t_1} \in \bigcup_{\substack{\Gf \in \C\left(\Gs\right)}}\Desc(C_{t_C},  \Gf\setminus \{ X^i_{t-\gamma_i}\})\setminus \Xf\right\}  =  \min \{t_1 \mid t_1 \geq t_C \text{ and } S_{t_1} \notin \Xf\}$.

    \item Else if $S \in \Desc(C, \Gs)$, then there exist directed paths from $C$ to $S$ in $\Gs$ but all of them goes through $X^i$. Thus for all $t_1 \geq t_C +1 $, there exists an FTCG $\Gf$  such that $\Gf \setminus \{ X^i_{t-\gamma_i} \}$ contains a directed path from $C_{t_C}$ to $S_{t_1}$. But there is no FTCG which contains a directed path from $C_{t_C}$ to $S_{t_C}$ because this path would need to go through $X^i_{t-\gamma_i}$. Therefore, $\min \left\{ t_1 \mid S_{t_1} \in \bigcup_{\substack{\Gf \in \C\left(\Gs\right)}}\Desc(C_{t_C}, \Gf\setminus \{ X^i_{t-\gamma_i}\})\setminus \Xf\right\}  =  \min \{t_1 \mid t_1 \geq t_C + 1 \text{ and } S_{t_1} \notin \Xf\}$.

    \item Otherwise, necessarily, for all $t_1$, $S_{t_1}$ cannot be a descendant of $C_{t_C}$ in any FTCG because otherwise S would be a descendant of $C$ in $\Gs$. Thus $\min \left\{ t_1 \mid S_{t_1} \in \bigcup_{\substack{\Gf \in \C\left(\Gs\right)}}\Desc(C_{t_C},  \Gf\setminus \{ X^i_{t-\gamma_i}\})\setminus \Xf\right\}  =  +\infty$.
\end{itemize}

Therefore, Equation \eqref{eq:B} holds.
\end{proof}

\begin{lemma}
\label{lemma:preuve_algo_t_NC}
Let $\Gs = (\mathcal{V}^s, \mathcal{E}^s)$ be an SCG and let $P(y_t \mid \text{do}(x^1_{t-\gamma_1}), \dots, \text{do}(x^n_{t-\gamma_n}))$ be the considered effect. Algorithm \ref{algo:calcul_t_NC} computes $(t_{\NC}(S))_{S \in \mathcal{V}^s}$ in pseudo-linear complexity. 
\end{lemma}

\begin{proof}
\textbf{Let us prove that Algorithm \ref{algo:calcul_t_NC} is correct:} \newline

Algorithm \ref{algo:calcul_t_NC} is split into two parts, the first one computes $t_C$ for all $C \in \Ch(\Xs)$ as defined in Lemma \ref{lemma:characterisation_CF}. The second part computes $(t_{\NC}(S))_{S \in \mathcal{V}^s}$. Lemma \ref{lemma:calcul_t_C} shows that Algorithm \ref{algo:calcul_t_NC} computes $(t_{\NC}(S))_{S \in \mathcal{V}^s}$ correctly. 

Let us prove that the second part of the algorithm computes  $(t_{\NC}(S))_{S \in \mathcal{V}^s}$ correctly. Let $S$ be a time series. By Lemma \ref{lemma:characterisation_CF}, we know that:
\[
    \NC = \bigcup_{\substack{C ~\in~ \Ch\left(\Xs\right)  \\ t_C~ <~ +\infty}}
        \ \ \bigcup_{\substack{\Gf \in \C\left(\Gs\right) \\ \exists i ~ \Gf \models X^i_{t-\gamma_i} \rightarrow C_{t_C}}}
        \Desc\left(C_{t_C}, \Gf \right) \setminus \Xf.
\]

Thus, we have:

\[
    t_{\NC}[S] =
        \min_{\substack{C ~\in~ \Ch\left(\Xs\right)  \\ t_C~ <~ +\infty}} 
            \min \left\{ 
                t_1 \mid S_{t_1} \in 
                    \bigcup_{\substack{\Gf \in \C\left(\Gs\right) \\ \exists i ~ \Gf \models X^i_{t-\gamma_i} \rightarrow C_{t_C}}}\Desc(C_{t_C}, \Gf)\setminus \Xf
            \right\}
\]

By Lemma \ref{lemma:characterisation_CF_suite} we have:

\[
t_{\NC}[S] =
        \min_{\substack{C ~\in~ \Ch\left(\Xs\right)  \\ t_C~ <~ +\infty}} 
            \min \left\{ 
                t_1 \mid S_{t_1} \in 
                    \begin{cases}
                        \bigcup_{\substack{\Gf \in \C\left(\Gs\right)}} \Desc\left(C_{t_C}, \Gf \right) & \text{if } d(C), \\
                        \bigcup_{\substack{\Gf \in \C\left(\Gs\right)}} \Desc\left(C_{t_C}, \Gf \setminus \{X^i_{t-\gamma_i}\} \right)  & \text{otherwise.}
                    \end{cases}
            \right\}
\]

Where $d(C) \coloneqq 
                     \# \{i \mid X^i \in \text{Pa}(C, \Gs) \text{ and }t-\gamma_i < t_C \} \geq 1 \text{ or }  
                     \# \{i \mid X^i \in \text{Pa}(C, \Gs) \text{ and } t-\gamma_i =t_C \} \geq 2$, and, in the second case, $X^i$ is the only element of $\{X^i \in \text{Pa}(C, \Gs) \text{ and } t-\gamma_i \leq t_C \} = \{X^i \in \text{Pa}(C, \Gs) \text{ and } t-\gamma_i = t_C \}$.
Thus, we have:

\begin{equation}
\label{eq:magie_primaire}
t_{\NC}[S] =
        \min_{\substack{C ~\in~ \Ch\left(\Xs\right)  \\ t_C~ <~ +\infty}} 
            \begin{cases}
                  \min \left\{ t_1 \mid S_{t_1} \in  \bigcup_{\substack{\Gf \in \C\left(\Gs\right)}} \Desc\left(C_{t_C}, \Gf \right) \right\}  & \text{if } d(C), \\
                \min \left\{ t_1 \mid S_{t_1} \in  \bigcup_{\substack{\Gf \in \C\left(\Gs\right)}} \Desc\left(C_{t_C}, \Gf \setminus \{X^i_{t-\gamma_i}\} \right) \right\} & \text{otherwise.}
            \end{cases}
\end{equation}

We can use Equation \eqref{eq:B} from Lemma \ref{lemma:truc} to specify the first case of Equation \eqref{eq:magie_primaire}. For the second case, we know that $X^i$ is the only element of $\{X^i \in \text{Pa}(C, \Gs) \text{ and } t-\gamma_i \leq t_C \} = \{X^i \in \text{Pa}(C, \Gs) \text{ and } t-\gamma_i = t_C \}$. Thus $t_C = t-\gamma_i$, and we can use Equation \eqref{eq:B} to specify this case. Algorithm \ref{algo:calcul_t_NC} uses these specifications to compute Equation \eqref{eq:magie_primaire}. It enumerates the $C \in \Ch(\Xs)$ such that $t_C < +\infty$ in the correct order to avoid unnecessary operations. Therefore, Algorithm \ref{algo:calcul_t_NC} is correct.

\textbf{We now prove that Algorithm \ref{algo:calcul_t_NC} works in pseudo linear complexity:} \newline
We will go through the lines of Algorithm \ref{algo:calcul_t_NC} and we will show that they run in pseudo-linear complexity:
\begin{itemize}
    \item The first line of Algorithm \ref{algo:calcul_t_NC} is $AncY \gets \left(\max \left\{ t_1 \mid  \exists\Gf\text{ s.t }  S_{t_1} \in \Anc(Y_{t}, \Gf \setminus \Xf)  \right\}\right)_{S \in \mathcal{V}^s}$. By Lemma \ref{lemma:calcul_t_NC_aux}, this line is computed by Algorithm \ref{algo:calcul_t_NC_aux} in $\mathcal{O} ((\left| \mathcal{E}^s \right| + \left| \mathcal{V}^s \right|) \log \left| \mathcal{V}^s \right|)$.

    \item Then, the algorithm computes $t_C$ for all $C \in \Ch(\Xs)$. While enumerating all the $C \in \Ch(\Xs)$, only the arrows from $\Xs$ are seen. Thus, this enumeration is done in $\mathcal{O} (\left| \mathcal{E}^s \right| )$ time.  With the appropriate pre-computations all the computations for each $C$ are done in $\mathcal{O}(1)$ time. Thus, this loop takes at most $\mathcal{O} (\left| \mathcal{E}^s \right| )$ time. 

    \item The sort of $L$ can be done in a complexity of $\mathcal{O} (\left| \mathcal{V}^s \right| \log \left| \mathcal{V}^s \right|)$. Indeed, $L$ contains at most $\left| \mathcal{V}^s \right|$ elements.

    \item Then, the algorithm computes $(t_{\NC}(S))_{S \in \mathcal{V}^s}$. For each $D$ computed by the algorithm, the algorithm computes $t_{\NC}(D) \gets \min \{t_1 \mid t_1 \geq t_\star \text{ and } D_{t_1} \notin \Xf\}$, where $t_\star \in \{t_C,t_C +1 \}$ depending on the context. These lines can be computed in $\mathcal{O}(1)$. Indeed, we need to distinguish two cases:
    \begin{itemize}
        \item If $D_{t_\star} \notin \Xf$, the computations becomes $t_{\NC}(D) \gets t_\star$. 
        \item Otherwise, we pick the value from $\{t_1 \mid t_1 \geq t-\gamma_i \text{ and } X^i_{t_1} \notin \Xf\}_{i \in \{1,\cdots,n\}}$, already pre-computed by Algorithm \ref{algo:calcul_t_NC:1}.
    \end{itemize}
    Since Algorithm \ref{algo:calcul_t_NC:1} runs in $\mathcal{O}(\left| \Xf \right| \log \left| \Xf \right|)$\footnote{The $\log n$ comes from sorting interventions by time index.}, all computations in all foreach loops are done in $\mathcal{O}(1)$ by adding $\mathcal{O}(\left| \Xf \right| \log \left| \Xf \right|)$ to the overall complexity of Algorithm \ref{algo:calcul_t_NC}. Moreover, by sharing the unseen set among all calculations of descendant sets, all time series that are not a child of some $X^i \in \Xs$ are seen at most one time and all $C \in \Ch(\Xs)$ are seen at most three times. Similarly, arrows from a time series that is not a child of some $X^i \in \Xs$ are seen at most one time and arrows from all $C \in \Ch(\Xs)$ are seen at most two times. Thus, the computation of $(t_{\NC}(S))_{S \in \mathcal{V}^s}$ is done in $\mathcal{O}(\left| \mathcal{E}^s \right| + \left| \mathcal{V}^s \right|)$.
\end{itemize}

Therefore,the overall complexity of Algorithm \ref{algo:calcul_t_NC} is $\mathcal{O}\left(\left| \Xf \right| \log \left| \Xf \right| + (\left| \mathcal{E}^s \right| + \left| \mathcal{V}^s \right|) \log \left| \mathcal{V}^s \right|\right)$. It is indeed a pseudo linear complexity with respect to $\Gs$ and $\Xf$.

\begin{algorithm}[H]
    \caption{Computation of $\{t_1 \mid t_1 \geq t-\gamma_i \text{ and } X^i_{t_1} \notin \Xf\}_{i \in \{1,\cdots,n\}}$}
    \label{algo:calcul_t_NC:1}
    \SetKwInput{KwData}{Input}
    \SetKwInput{KwResult}{Output}
    
    \KwData{List of interventions, sorted by decreasing time indices.}
    \KwResult{$\{ \min \{t_1 \mid t_1 \geq t-\gamma_i \text{ and } X^i_{t_1} \notin \Xf\}\}_{i \in \{1,\cdots,n\}}$, denoted as $\{t_{X^i_{t-\gamma_i}}\}_i$}
    
    $L \gets$ List of lists of interventions, grouped by time series, preserving the time index ordering\;
    
    \For{$l \in L$}{
        $X^i_{t-\gamma_i} \gets l[0]$\;
        $t_{X^i_{t-\gamma_i}} \gets t - \gamma_i + 1$\;
        
        \ForEach{$X^j_{t-\gamma_j}$ in $l[1:]$}{
            Let $X^i_{t-\gamma_i}$ be the predecessor of $X^j_{t-\gamma_j}$ in $l$\;
            \eIf{$t-\gamma_i +1 = t-\gamma_j$}{
                $t_{X^j_{t-\gamma_j}} \gets t_{X^i_{t-\gamma_i}}$\;
            }{
                $t_{X^j_{t-\gamma_j}} \gets t-\gamma_j + 1$\;
            }
        }
    }
\end{algorithm}
\end{proof}

\begin{lemma}
\label{lemma:calcul_t_NC_aux}
    Let $\Gs = (\mathcal{V}^s, \mathcal{E}^s)$ be an SCG and let $P(y_t \mid \text{do}(x^1_{t-\gamma_1}), \dots, \text{do}(x^n_{t-\gamma_n}))$ be the considered effect. $\left(\max \left\{ t_1 \mid  \exists\Gf\text{ s.t } 
      S_{t_1} \in \Anc(Y_{t}, \Gf \setminus \Xf)  \right\}\right)_{S \in \mathcal{V}^s}$ is computed by Algorithm \ref{algo:calcul_t_NC_aux} in $\mathcal{O} ((\left| \mathcal{E}^s \right| + \left| \mathcal{V}^s \right|) \log \left| \mathcal{V}^s \right|)$.
\end{lemma}

\begin{algorithm}[t]
\setlength{\rightskip}{-.5cm} 
\caption{Computation of $\left(\max \left\{ t_1 \mid  \exists\Gf\text{ s.t } 
  S_{t_1} \in \Anc(Y_{t}, \Gf \setminus \Xf)  \right\}\right)_{S \in \mathcal{V}^s}$}\label{algo:calcul_t_NC_aux}
\SetKwInput{KwData}{Input}
\SetKwInput{KwResult}{Output}
\KwData{ $\Gs = (\mathcal{V}^s, \mathcal{E}^s)$ an SCG, $\Xf$ and $V_{t_v} \in \mathcal{V}^f$}
\KwResult{$\left(\max \left\{ t_1 \mid  \exists\Gf\text{ s.t } 
  S_{t_1} \in \Anc(Y_{t}, \Gf \setminus \Xf)  \right\}\right)_{S \in \mathcal{V}^s}$ written as $(AncY[S])_{S \in \mathcal{V}^s}.$}
$Q \gets \text{PriorityQueue}(Y_t)$ \;
$AncY[S] \gets -\infty \quad \forall S \in \mathcal{V}^s$\;
$S.seen \gets False \quad \forall S \in \mathcal{V}^s$\;
$AncY[Y] \gets 0$\;
$Y.seen \gets True$
\While{$Q \neq \emptyset$}{
    $S_{t_s} \gets Q.\text{pop\_element\_with\_max\_time\_index}()$\;
    \ForEach{unseen $P \in \text{Pa}(S, \Gs)$}{
        $AncY[P] \gets \max \{t_1 \mid t_1 \leq t_s \text{ and } P_{t_1} \notin \Xf\}$\;
        \lIf{$AncY[P] \neq -\infty$}{
            $Q.insert(P_{AncY[P]})$
        }
        $P.seen \gets true$ \;
    }
}
\end{algorithm}

\begin{proof}
    The proof is identical to the proof of Lemma \ref{lemma:calcul_V_E_acc} except that we compute the $Y_t$-$(\mathcal{V}^s \setminus \Xf)$-accessibility instead of the $Y_t$-$(\CF \setminus \Xf)$-accessibility. Precisely, in this case, $Y_t$ is considered to be $Y_t$-$\mathcal{V}^s \setminus \Xf$-accessible.
\end{proof}

\subsubsection{Proof of Theorem \ref{th:equiv_IBC_multivarie}}

\mytheoremequivIBC*

\begin{proof}
    \textbf{We start by proving that every proper non-causal paths from $\Xf$ to $Y_t$ that leaves $\CF$, except perhaps for its first vertex, are blocked by $\C \coloneqq \left(\mathcal{V}^f \setminus \NC \right) \setminus X^f$:} 
    
    Let us consider an FTCG $\Gf$ and a proper non-causal path $\pi^f$ from $X^i_{t-\gamma_i} \in \Xf$ to $Y_t$ that leaves $\CF \cup \{ X^i_{t-\gamma_i}\}$. Let $C_{t_c}$ be the last vertex of $\pi^f$ in $\mathcal{V}^f \setminus \left(\CF \cup \{ X^i_{t-\gamma_i}\}\right)$. We know that $Y$ is a descendant of an intervention in $\Gs$. Thus, $Y_t \in \CF$. Thus, $C_{t_c} \neq Y_t$ and $C_{t_C}$ is not the last vertex of $\pi^f$. Let us consider $D_{t_d}$, the successor of $C_{t_c}$ in $\pi^f$. By definition, $D_{t_d} \in \CF$. Necessarily $\pi^f_{\mid \{C_{t_c}, D_{t_d}\}}  = D_{t_d} \leftarrow C_{t_c}$. Indeed, otherwise $C_{t_c}$ would be in $\CF$. Hence, $C_{t_c}$ is not a collider on $\pi^f$. Moreover, $\pi^f$ is a proper path thus, $C_{t_c} \in  \left(\mathcal{V}^f \setminus \CF \right) \setminus X^f =  \left(\mathcal{V}^f \setminus \NC \right) \setminus X^f$. Therefore $\C = \left(\mathcal{V}^f \setminus \NC \right) \setminus X^f$ blocks $\pi^f$.

    \noindent \textbf{Let us prove the two implications of the theorem:}
    \begin{itemize}
        \item $\ref{th:equiv_IBC_multivarie:1} \Rightarrow \ref{th:equiv_IBC_multivarie:2}$: Let us prove the contrapositive: Let $\Gf$ be an FTCG that contains $\pi^f$, a collider-free backdoor path from $X^i_{t-\gamma_i}$ to $Y_t$ that remains within $\NC \cup \{ X^i_{t-\gamma_i}\}$. Thus $\pi^f$ is a proper non-causal path from $\Xf$ to $Y_t$. Moreover, $\pi^f$ does not have a collider, so the only way to block it is by conditioning on one of its vertices. However, all its vertices are within $\NC \cup \{ X^i_{t-\gamma_i}\} \subseteq \NC \cup \Xf$. Thus, $\pi^f$ cannot be blocked. Therefore, the effect is not identifiable by common adjustment.
        
        \item $\ref{th:equiv_IBC_multivarie:2} \Rightarrow \ref{th:equiv_IBC_multivarie:1}$: Let us assume that condition \ref{th:equiv_IBC_multivarie:2} holds. We will show that $\C = \left(\mathcal{V}^f \setminus \NC \right) \setminus X^f$ is a valid common adjustment set. Let $\Gf \in \C(\Gs)$ be an FTCG. Let us check that $\C$ is a valid adjustment set in $\Gf$:
        \begin{itemize}
            \item By definition, $\C \cap \CF = \emptyset$. Thus $\C \cap \Forb(\Xf, Y_t, \Gf) = \emptyset$.
            \item Let $\pi^f$ be a proper non-causal path from $\Xf$ to $Y_t$ in $\Gf$. Let $X^i_{t-\gamma_i}$ be the first vertex of $\pi^f$. We distinguish two cases:
            \begin{itemize}
                \item If $\pi^f$ leaves $\CF \cup \{X^i_{t-\gamma_i} \}$, then $\C$ blocks $\pi^f$.
                \item Otherwise, $\pi^f$ remains within $\CF \cup \{X^i_{t-\gamma_i} \}$. $\pi^f$ is a proper path, hence $\pi^f$ remains within  $\NC \cup \{X^i_{t-\gamma_i} \}$. By contradiction, we show that $\pi^f$ contains a collider. Let us assume that $\pi^f$ is collider-free. We distinguish two cases:
                \begin{itemize}
                    \item If $\pi^f$ starts by $X^i_{t-\gamma_i} \rightarrow$, then $\pi^f$ would be a causal path, which contradicts the definition of $\pi^f$.
                    \item If $\pi^f$ starts by $X^i_{t-\gamma_i} \leftarrow$, then $\pi^f$ would be a collider-free backdoor path from $X^i_{t-\gamma_i}$ to $Y_t$ that remains in $\NC \cup \{ X^i_{t-\gamma_i}\}$, which contradicts condition \ref{th:equiv_IBC_multivarie:2}.
                \end{itemize}
                Therefore, $\pi^f$ contains a collider denoted $C_{t_c}$. $C_{t_c} \in \CF$, thus every descendant of $C_{t_c}$ belongs to $\CF$. Thus $\C \cap \Desc(C_{t_c}, \Gf) \subseteq \C \cap \CF = \emptyset$. Thus $\C$ blocks $\pi^f$.
            \end{itemize}
            In all cases, $\C$ blocks $\pi^f$.
        \end{itemize}
        Therefore, $\C$ is a valid adjustment set in $\Gf$. Therefore, $\C = \left(\mathcal{V}^f \setminus \NC \right) \setminus X^f$ is a valid common adjustment set.
    \end{itemize}
\end{proof}

\subsection{Proofs of Section \ref{ssct:in practice}}
\subsubsection{Proofs of Section \ref{subsubsection:CfbWF}}

\mylemmaIBCenumerationcheminsdiriges*

\begin{proof}
    \noindent \textbf{Let us prove the two implications:}
    \begin{itemize}
        \item $\ref{lemma:IBC_enumeration_chemins_diriges:1} \Rightarrow \ref{lemma:IBC_enumeration_chemins_diriges:2}$: Let us consider $X^i_{t-\gamma_i}$ and an FTCG $\Gf$ belonging to $\C(\Gs)$  such that $\Gf$ contains  $ \pi^f \coloneqq X^i_{t- \gamma_i} \leftsquigarrow Y_t$ which remains in $\NC \cup \{ X^i_{t-\gamma_i}\}$. By definition, we already know that $t- \gamma_i \leq t$. Causality does not move backwards in time thus $t- \gamma_i \geq t$. Therefore, $t-\gamma_i = t$ and all vertices of $\pi^f$ are a time $t$. Therefore the reduction $\pi^s$ of $\pi^f$ is still a path in $\Gs$. Moreover, $\pi^f$ remains in $\NC \cup \{X^i_{t-\gamma_i}\}$ thus $\pi^s \subseteq \mathcal{S}$. Therefore, $X^i_{t-\gamma_i} \in \Desc(Y, \Gs_{\mid \mathcal{S}})$ and the intervention is at time t.
        
        \item $\ref{lemma:IBC_enumeration_chemins_diriges:2} \Rightarrow \ref{lemma:IBC_enumeration_chemins_diriges:1}$: Let us consider $\pi^s = X^i \leftarrow V^2 \leftarrow \cdots \leftarrow V^{n-1} \leftarrow Y$, the smallest path from $Y$ to an interventional time series at time $t$ in $\Gs_{\mid \mathcal{S}}$. For all $i \in \{2, \cdots, n-1 \}$, $V^i_t$ is not an intervention because otherwise we could find a smaller path that $\pi^f$. We can construct an FTCG $\Gf$ which contains the path $\pi^f = X^i_t \leftarrow V^2_t \leftarrow \cdots \leftarrow V^{n-1}_t \leftarrow Y_t$. $\pi^f$ remains in $\NC \cup \{X^i_{t-\gamma_i}\}$ because for all $i$, $t_{\NC}(V^i) \leq 0$. Therefore, there exist an intervention $X^i_{t-\gamma_i}$ and an FTCG $\Gf$ belonging to $\C(\Gs)$ which contains $X^i_{t- \gamma_i} \leftsquigarrow Y_t$ which remains in $\NC \cup \{X^i_{t-\gamma_i}\}$.
    \end{itemize}
\end{proof}

\subsubsection{Proofs of Section \ref{subsubsection:CfbF}}

\mylemmaequivexistencecheminforkCDwithoutconsistencythroughtime*

\begin{proof}
    \noindent \textbf{Let us prove the two implications of the theorem:}
    \begin{itemize}
        \item $\ref{lemma:equiv_existence_chemin_fork_CD_without_consistency_through_time:1}
        \Rightarrow \ref{lemma:equiv_existence_chemin_fork_CD_without_consistency_through_time:2}$:
        Let $X^i_{t-\gamma_i}$  be the intervention and $\Gf$ the FTCG belonging to $\C(\Gs)$ which contains the path $X^i_{t- \gamma_i} \leftsquigarrow F_{t_f} \rightsquigarrow Y_t$ which remains in $\NC \cup \{X^i_{t-\gamma_i}\}$. $\Gf$ proves that $F_{t_f}$ is $X^i_{t- \gamma_i}$-$\NC$-accessible and $Y_t$-$\NC$-accessible.
        \item $\ref{lemma:equiv_existence_chemin_fork_CD_without_consistency_through_time:2}
        \Rightarrow \ref{lemma:equiv_existence_chemin_fork_CD_without_consistency_through_time:1}$:
        Let $X^i_{t- \gamma_i}$ and $F_{t_f}$ be such that $F_{t_f}$ is $X^i_{t- \gamma_i}$-$\NC$-accessible and $Y_t$-$\NC$-accessible. Therefore there exists $\Gf_1$ in which there is $\pi^f_1 \vcentcolon= X^i_{t- \gamma_i} \leftsquigarrow F_{t_f}$ which remains in $\NC \cup \{X^i_{t-\gamma_i}\}$and there exists $\Gf_2$ in which there is $\pi^f_2 \vcentcolon= F_{t_f} \rightsquigarrow Y_t$ which remains in $\NC$. We distinguish two cases:
        \begin{itemize}
            \item If $\pi^f_1 \cap \pi^f_2 = \{ F_{t_f}\}$, then we can build an FTCG $\Gf_3 \in \C(\Gs)$ which contains $\pi^f_3 \coloneqq X^i_{t- \gamma_i} \leftsquigarrow F_{t_f} \rightsquigarrow Y_t$ the concatenation of $\pi^f_1$ and $\pi^f_2$. $\pi^f_3$ is a path without cycle because $\pi^f_1 \cap \pi^f_2 = \{ F_{t_f}\}$. It remains in $\NC \cup \{X^i_{t-\gamma_i}\}$ because its vertices come from $\pi^f_1$ and $\pi^f_2$ which remain in $\NC\cup \{X^i_{t-\gamma_i}\}$.
            \item Otherwise, $\pi^f_1 \cap \pi^f_2$ contains at least two element. We know that $X^i_{t-\gamma_i} \notin \pi^f_1 \cap \pi^f_2$ because $\pi^f_2$ remains in $\NC$. Similarly, we know that $Y_t \notin \pi^f_1 \cap \pi^f_2$ because there would be a directed path from $Y_t$ to $X^i_{t-\gamma_i}$. Let us consider $V_{t_v}$ the latest element of $\pi^f_1$ in $\pi^f_1 \cap \pi^f_2$. In $\Gf_1$ there is ${\pi^f_1}' \coloneqq X^i_{t- \gamma_i} \leftsquigarrow V_{t_v}$, in $\Gf_2$ there is ${\pi^f_2}' \coloneqq  V_{t_v} \rightsquigarrow  Y_t$ and ${\pi^f_1}' \cap {\pi^f_2}' = \{ V_{t_v}\}$.Therefore, with the same reasoning as in the first case, we can construct $\Gf_3$ which contains the path $X^i_{t - \gamma_i} \leftsquigarrow V_{t_v} \rightsquigarrow Y_t$ which remains in $\NC \cup \{X^i_{t-\gamma_i}\}$.
        \end{itemize}
        In all cases, there exist an intervention $X^i_{t-\gamma_i}$ and an FTCG $\Gf$ belonging to $\C(\Gs)$ which contains the path $X^i_{t- \gamma_i} \leftsquigarrow F_{t_f} \rightsquigarrow Y_t$ which remains in $\NC \cup \{X^i_{t-\gamma_i}\}$.
    \end{itemize}
\end{proof}

As shown by Lemma~\ref{lemma:CharactVNCAcc}, for any time series \(F\in\mathcal{V}^s\), knowledge of $t^{\NC}_{V_{t_v}}(F)$ and $t_{\NC}(F)$ is sufficient to characterize efficiently the set $\{ t_1 \mid F_{t_1} \text{is $V_{t_v}$-$\NC$-accessible} \}$.

\begin{restatable}{lemma}{mylemmaCharactVNCAcc}{}
\label{lemma:CharactVNCAcc}
    Let $\Gs = (\mathcal{V}^s, \mathcal{E}^s)$ be an SCG, $P(y_t \mid \text{do}(x^1_{t-\gamma_1}), \dots, \Do(x^n_{t-\gamma_n}))$ be the considered effect, $V_{t_v}\in \mathcal{V}^f$  and $F\in \mathcal{V}^s$. Then the following statements are equivalent:
    \begin{enumerate}
        \item $F_{t_f}$ is $V_{t_v}$-$\NC$-accessible.
        \label{lemma:CharactVNCAcc:1}
        \item $t_{\NC}(F) \leq t_f \leq t^{\NC}_{V_{t_v}}(F)$ and $F_{t_f} \notin \Xf$.
        \label{lemma:CharactVNCAcc:2}
    \end{enumerate}
\end{restatable}

\begin{proof}
    \noindent \textbf{Let us prove the two implications:}
    \begin{itemize}
        \item $\ref{lemma:CharactVNCAcc:1} \Rightarrow \ref{lemma:CharactVNCAcc:2}$: Let $F_{t_f}$ and $V_{t_v}$ be such that $F_{t_f}$ is $V_{t_v}$-$\NC$-accessible. Thus $F_{t_f} \in \NC$. Hence, $F_{t_f} \notin \Xf$ and $t_{\NC}(F) \leq t_f$ by definition of $t_{\NC}(F)$. By definition of $t^{\NC}_{V_{t_v}}(F)$, $ t_f \leq t^{\NC}_{V_{t_v}}(F)$. Therefore $t_{\NC}(F) \leq t_f \leq t^{\NC}_{V_{t_v}}(F)$ and $F_{t_f} \notin \Xf$.

        \item $\ref{lemma:CharactVNCAcc:1} \Rightarrow \ref{lemma:CharactVNCAcc:2}$: Let $F_{t_f}$ and $V_{t_v}$ be such that $t_{\NC}(F) \leq t_f \leq t^{\NC}_{V_{t_v}}(F)$ and $F_{t_f} \notin \{ X^i_{t-\gamma_i} \}_i$. Since $t_{\NC}(F) \leq t_f$, it follows that $t_{\NC}(F) \neq +\infty$. By Lemma \ref{lemma:def_equiv_CD}, we have $\NC = \bigcup_{Z \in \mathcal{V}^S} \{Z_{t_1}\}_{t_1 \geq t_{\NC}(Z)} \setminus \{ X^i_{t-\gamma_i} \}_i$, and given that $F_{t_f} \notin \{ X^i_{t-\gamma_i} \}_i$, it follows that $F_{t_f} \in \NC$. Since $t_f \leq t^{\NC}_{V_{t_v}}(F)$, it follows that $t^{\NC}_{V_{t_v}}(F) \neq -\infty$. Thus, $F_{t^{\NC}_{V_{t_v}}(F)}$ is $V_{t_v}$-$\NC$-accessible. Hence, there exists an FTCG in which there is a path $ \pi^f: V_{t_v} \leftsquigarrow F_{t^{\NC}_{V_{t_v}}(F)}$ that remains in $\NC$. By changing the last arrow of this path, we can construct an FTCG which contains the path $V_{t_v} \leftsquigarrow F_{t_f}$. This path remains in $\NC$ because $F_{t_f} \in \NC$ and all other vertices of the paths are in $\pi^f$ which remains in $\NC$. Therefore $F_{t_f}$ is $V_{t_v}$-$\NC$-accessible.

    \end{itemize}
\end{proof}

By setting \(t_f = t_{\NC}(F)\) in Lemma~\ref{lemma:CharactVNCAcc}, we observe that it suffices to test only the single time point \(F_{t_{\NC}(F)}\) in Condition 2 of Lemma~\ref{lemma:equiv_existence_chemin_fork_CD_without_consistency_through_time}.

\begin{restatable}{corollary}{mycorUn}{}
\label{cor:1}
    Let $\Gs = (\mathcal{V}^s, \mathcal{E}^s)$ be an SCG,  $P(y_t \mid \text{do}(x^1_{t-\gamma_1}), \dots, \Do(x^n_{t-\gamma_n}))$ be the considered effect, $V_{t_v} \in \mathcal{V}^f$ and $F \in \mathcal{V}^s$. Let $X^i_{t - \gamma_i}$ be a fixed intervention. The following statements are equivalent:
    \begin{enumerate}
        \item There exists $t_f$ such that $F_{t_f}$ is $X^i_{t - \gamma_i}$-$\NC$-accessible and $Y_t$-$\NC$-accessible.
        \label{cor:1:1}
        \item $t_{\NC}(F) \leq t^{\NC}_{X^i_{t - \gamma_i}}(F)$ and $t_{\NC}(F) \leq t^{\NC}_{Y_t}(F)$.
        \label{cor:1:2}
    \end{enumerate}
\end{restatable}

\begin{proof}
    Let us show the two implications of the corollary:
    \begin{itemize}
        \item $\ref{cor:1:1} \Rightarrow \ref{cor:1:2}$: If there exists $t_f$ such that $F_{t_f}$ is $X^i_{t - \gamma_i}$-$\NC$-accessible and $Y_t$-$\NC$-accessible then $F_{t_{\NC}(F)}$ is $X^i_{t - \gamma_i}$-$\NC$-accessible and $Y_t$-$\NC$-accessible. Thus, by Lemma \ref{lemma:CharactVNCAcc}, $t_{\NC}(F) \leq t_{\NC}(F) \leq t^{\NC}_{X^i_{t - \gamma_i}}(F)$ and $F_{t_{\NC}(F)} \notin \Xf$ and $t_{\NC}(F) \leq t_{\NC}(F) \leq t^{\NC}_{Y_t}(F)$ and $F_{t_{\NC}(F)} \notin \Xf$. Hence $t_{\NC}(F) \leq t^{\NC}_{X^i_{t - \gamma_i}}(F)$ and $t_{\NC}(F) \leq t^{\NC}_{Y_t}(F)$.

        \item $\ref{cor:1:2} \Rightarrow \ref{cor:1:1}$: If $t_{\NC}(F) \leq t^{\NC}_{X^i_{t - \gamma_i}}(F)$ and $t_{\NC}(F) \leq t^{\NC}_{Y_t}(F)$, then $t_{\NC}(F) \neq + \infty$ and $t^{\NC}_{X^i_{t - \gamma_i}}(F) \neq -\infty$. Therefore, $F_{t_{\NC}(F)} \notin \Xf$ and  $t_{\NC}(F) \leq t_{\NC}(F) \leq t^{\NC}_{X^i_{t - \gamma_i}}(F)$ and $F_{t_{\NC}(F)} \notin \Xf$ and $t_{\NC}(F) \leq t_{\NC}(F) \leq t^{\NC}_{Y_t}(F)$. Therefore, by Lemma \ref{lemma:CharactVNCAcc}, $F_{t_{\NC}(F)}$ is $X^i_{t - \gamma_i}$-$\NC$-accessible and $Y_t$-$\NC$-accessible. Therefore, there exists $t_f$ such that $F_{t_f}$ is $X^i_{t - \gamma_i}$-$\NC$-accessible and $Y_t$-$\NC$-accessible.
    \end{itemize}
    
\end{proof}

\subsubsection{Proofs of Section \ref{sssct:algo_IBC}}

\begin{restatable}{lemma}{mylemmacalculVEacc}{}
    \label{lemma:calcul_V_E_acc}
    Let $\Gs = (\mathcal{V}^s, \mathcal{E}^s)$ be an SCG,  $P(y_t \mid \text{do}(x^1_{t-\gamma_1}), \dots, \Do(x^n_{t-\gamma_n}))$ be the considered effect and $V_{t_v}\in \mathcal{V}^f$. The set $\{t^{\NC}_{V_{t_v}}(F) \mid F \in \mathcal{V}^S \}$ {can be} computed by Algorithm \ref{algo:calcul_V_E_acc} {which} complexity is $\mathcal{O} \left(\left| \mathcal{E}^s \right| + \left| \mathcal{V}^s \right| \log \left| \mathcal{V}^s \right| \right)$\footnote{In amortized time. Using a binary heap, the algorithm runs in $\mathcal{O} ((\left| \mathcal{E}^s \right| + \left| \mathcal{V}^s \right|) \log \left| \mathcal{V}^s \right|)$.}.
\end{restatable}

\begin{proof}
    \noindent Let $V_{t_v}$ be a temporal variable.\newline
    
    \textbf{Let us prove that Algorithm \ref{algo:calcul_V_E_acc} terminates:}

    \noindent Algorithm \ref{algo:calcul_V_E_acc} \textbf{terminates} because at each step of the \textbf{while} loop, (number of unseen times series, length of $Q$) is strictly decreasing with respect to the {lexicographic ordering}.\newline

    \textbf{Let us prove that Algorithm \ref{algo:calcul_V_E_acc} is correct:}

             Firstly, by induction, \textbf{we show that the algorithm computes the correct value for each seen time series:}
        \begin{itemize}
            \item At the first step of the while loop, for each parent $P$ of $V$ in $\Gs$, the algorithm will compute $t^{\NC}_{V_{t_v}}(P) \gets \max \{t_1 \mid t_1 \leq t_v \text{ and } P_{t_1} \in \NC \setminus \{V_{t_v}\}\}$.\newline
            
            Let $P$ be a parent of $V$ in $\Gs$. We will show that $\{ t_1 \mid P_{t_1} \text{is $V_{t_v}$-$\NC$-accessible} \}  = \{t_1 \mid t_1 \leq t_v \text{ and } P_{t_1} \in \NC \setminus \{V_{t_v}\}\}$ by showing the two inclusions:
            \begin{itemize}
                \item For all $t_1 \in \{t_1 \mid t_1 \leq t_v \text{ and } P_{t_1} \in \NC \setminus \{V_{t_v}\}\}$, we can construct an FTCG $\Gf_{t_1} \in \C(\Gs)$ which contains $P_{t_1} \rightarrow V_{t_v}$. Therefore, all $P_{t_1} \in \{ P_{t_1} \mid t_1 \leq t_v \text{ and } P_{t_1} \in \NC \setminus \{V_{t_v}\}\}$ are $V_{t_v}$-$\NC$-accessible, \textit{i.e.}, $\{t_1 \mid t_1 \leq t_v \text{ and } P_{t_1} \in \NC \setminus \{V_{t_v}\}\} \subseteq \{ t_1 \mid P_{t_1} \text{is $V_{t_v}$-$\NC$-accessible} \}$

                \item  Let $t_1$ be such that $P_{t_1}$ is $V_{t_v}$-$\NC$-accessible. Causality does not move backwards in time thus $t_1 \leq t_v$. Moreover, $P_{t_1} \in \NC \setminus \{V_{t_v}\}$. Indeed, $P_{t_1}$ is $V_{t_v}$-$\NC$-accessible thus $P_{t_1} \in \NC$ and $P_{t_1} \neq V_{t_v}$ because no FTCG can contain a self loop. Therefore $ \{ t_1 \mid P_{t_1} \text{is $V_{t_v}$-$\NC$-accessible} \} \subseteq \{t_1 \mid t_1 \leq t_v \text{ and } P_{t_1} \in \NC \setminus \{V_{t_v}\}\}$ 
            \end{itemize}
            Therefore $\{ t_1 \mid F_{t_1} \text{is $V_{t_v}$-$\NC$-accessible} \}  = \{t_1 \mid t_1 \leq t_v \text{ and } P_{t_1} \in \NC \setminus \{V_{t_v}\}\}$, thus $t^{\NC}_{V_{t_v}}(P) =\max \{t_1 \mid t_1 \leq t_v \text{ and } P_{t_1} \in \NC \setminus \{V_{t_v}\}\}$.\newline
            
            \textbf{Therefore the algorithm computes correct values at the first step of the loop.}

            \item Let us suppose that the algorithm is correct until the $(n-1)$-th loop step. It pops $S_{t_s} = S_{t^{\NC}_{V_{t_v}}(S)}$ with $t^{\NC}_{V_{t_v}}(S)$ being the element of $Q$ with maximum time index. Let $P$ be an unseen parent of $S$ in $\Gs$, the algorithm will compute $t^{\NC}_{V_{t_v}}(P) \gets \max \{t_1 \mid t_1 \leq t_s \text{ and } P_{t_1} \in \NC \setminus \{V_{t_v}\}\}$.  We will show that $\{ t_1 \mid P_{t_1} \text{is $V_{t_v}$-$\NC$-accessible} \}  = \{t_1 \mid t_1 \leq t_v \text{ and } P_{t_1} \in \NC \setminus \{V_{t_v}\}\}$ by showing the two inclusions:
            \begin{itemize}
                \item $S_{t^{\NC}_{V_{t_v}}(S)}$ is $\NC$-accessible thus there exists an FTCG $\Gf \in \C(\Gs)$ which contains $S_{t^{\NC}_{V_{t_v}}(S)} \rightsquigarrow V_{t_v}$. For all $t_1 \in \{t_1 \mid t_1 \leq t_s \text{ and } P_{t_1} \in \NC \setminus \{V_{t_v}\}\}$, we can construct an FTCG $\Gf_{t_1} \in \C(\Gs)$ which contains $P_{t_1} \rightarrow S_{t^{\NC}_{V_{t_v}}(S)} \rightsquigarrow V_{t_v}$. Therefore, all $P_{t_1} \in \{ P_{t_1} \mid t_1 \leq t_s \text{ and } P_{t_1} \in \NC \setminus \{V_{t_v}\}\}$ are $V_{t_v}$-$\NC$-accessible, \textit{i.e.}, $\{t_1 \mid t_1 \leq t_v \text{ and } P_{t_1} \in \NC \setminus \{V_{t_v}\}\} \subseteq \{ t_1 \mid P_{t_1} \text{is $V_{t_v}$-$\NC$-accessible} \}$

                \item  Let $t_1$ be such that $P_{t_1}$ is $V_{t_v}$-$\NC$-accessible. Thus, $P_{t_1} \in \NC \setminus \{V_{t_v}\}$. Let us show that $t_1 \leq t_s$: there exists an FTCG $\Gf$ which contains $\pi^f \coloneqq P_{t_1} \rightsquigarrow V_{t_v}$ which remains in $\NC$. Let $U_{t_u}$ be the successor of $P_{t_1}$ in $\pi^f$. Causality does not move backwards in time thus $t_1 \leq t_u$. We distinguish two cases:
                \begin{itemize}
                    \item If $U = S$, then we already have $t_1 \leq t_s$.

                    \item Otherwise, $t_u \leq t_s$. Indeed, otherwise, thanks to the priority queue, U would have been seen before $S$ and $P$ would not be an unseen vertex. Therefore $t_1 \leq t_u \leq t_s$.
                \end{itemize}
                In all cases $t_1 \leq t_s$. Therefore, $ \{ t_1 \mid P_{t_1} \text{is $V_{t_v}$-$\NC$-accessible} \} \subseteq \{t_1 \mid t_1 \leq t_s \text{ and } P_{t_1} \in \NC \setminus \{V_{t_v}\}\}$.
            \end{itemize}
            Therefore $\{ t_1 \mid P_{t_1} \text{is $V_{t_v}$-$\NC$-accessible} \}  = \{t_1 \mid t_1 \leq t_s \text{ and } P_{t_1} \in \NC \setminus \{V_{t_v}\}\}$, thus $t^{\NC}_{V_{t_v}}(P) =\max \{t_1 \mid t_1 \leq t_v \text{ and } P_{t_1} \in \NC \setminus \{V_{t_v}\}\}$.\newline
            
            \textbf{Therefore the algorithm computes correct values at the $n$-th step of the loop.}
        \end{itemize}
        \textbf{Therefore, by induction principle, the algorithm computes the correct value for each seen time series.}

        Secondly, \textbf{we show that the algorithm sees all time series $S$ such that $t^{\NC}_{V_{t_v}}(S) \neq -\infty$:}

        Let $S$ be a time series such that $t^{\NC}_{V_{t_v}}(S) \neq -\infty$. By definition, $S_{t^{\NC}_{V_{t_v}}(S)}$ is $V_{t_v}$-$\NC$-accessible. Therefore there exists an FTCG in which there exists a path $\pi^f \coloneqq S_{t^{\NC}_{V_{t_v}}(S)} \rightsquigarrow V_{t_v}$ which remains in $\NC$. By cutting unnecessary parts and using Lemma \ref{lemma:helper:construct_FTCG_with_path} or \ref{lemma:helper:construct_FTCG_with_path_with_cst}, we can assume that $\pi^f$ does not pass twice by the same time series except perhaps for $V$ if $S = V$. Thus there is a directed path $\pi^S$ from $S$ to $V$ in $\Gs$ in which all vertices $U$ have a finite $t^{\NC}_{V_{t_v}}(U)$. Since, if the algorithm sees $U$ and $t^{\NC}_{V_{t_v}}(U) \neq -\infty$ then all parents of $U$ are seen by the algorithm and since all the parents of $V$ are seen by the algorithm, by induction on $\pi^s$, we can conclude that $S$ is seen.

        Therefore, \textbf{the algorithm sees all time series $S$ such that $t^{\NC}_{V_{t_v}}(S) \neq -\infty$}.

        Finally, \textbf{we show that Algorithm \ref{algo:calcul_V_E_acc} is correct:}
        Let $S$ be a time series. We distinguish two cases:
        \begin{itemize}
            \item If $S$ is seen by the algorithm, then the correct value is computed by the algorithm.
            \item Otherwise, the algorithm computes $t^{\NC}_{V_{t_v}}(S) \gets - \infty$ thanks to the initialisation step, which is the correct value.
        \end{itemize}

        \textbf{Therefore, Algorithm \ref{algo:calcul_V_E_acc} is correct.}\newline
    
    \noindent \textbf{Let us prove that Algorithm \ref{algo:calcul_V_E_acc} has a complexity of  $\mathcal{O} \left(\left| \mathcal{E}^s \right| + \left| \mathcal{V}^s \right| \cdot \log  \left| \mathcal{V}^s \right|\right)$:}

    $t^{\NC}_{V_{t_v}}(P) \gets \max \{t_1 \mid t_1 \leq t_s \text{ and } P_{t_1} \in \NC \setminus \{V_{t_v}\}\}$ can be computed in $\mathcal{O}(1)$ during the algorithm. Indeed, to do so we need to run first an algorithm similar to Algorithm \ref{algo:calcul_t_NC:1}.

    Let us assume that the priority queue is implemented using a Fibonacci heap. In this case, inserting an element takes $\mathcal{O}(1)$ amortized time, and extracting the max element takes $\mathcal{O}(\log \mid Q \mid)$ amortized time. During the execution of the algorithm, the priority queue $Q$ contains at most $\left| \mathcal{V}^s \right|$ elements, where $\left| \mathcal{V}^s \right|$ is the number of time series (or nodes) in $\Gs$. Therefore, each extraction costs $\mathcal{O}(\log \left| \mathcal{V}^s \right|)$, and since there are at most $\left| \mathcal{V}^s \right|$ extractions, the total cost of all extractions is $\mathcal{O}(\left| \mathcal{V}^s \right| \log \left| \mathcal{V}^s \right|)$.
    Additionally, each edge in $\Gs$ is processed at most once during the computation of the unseen parents of $S$ in $\Gs$ and the update of $t^{\NC}{V{t_v}}(P)$. These computations are performed in $\mathcal{O}(1)$ time. Thus, the total cost of processing all edges is $\mathcal{O}(\left| \mathcal{E}^s \right|)$.

    \textbf{Therefore, the overall complexity of the algorithm is $\mathcal{O}(\left| \mathcal{E}^s \right| + \left| \mathcal{V}^s \right| \log \left| \mathcal{V}^s \right|)$.}
    
\end{proof}

\mythforalgoIBC*

\begin{proof}
    Algorithm \ref{algo:calcul_IBC} uses directly the characterizations of Lemma \ref{lemma:IBC_enumeration_chemins_diriges} and Lemma \ref{lemma:equiv_existence_chemin_fork_CD_without_consistency_through_time} via Corollary~\ref{cor:1}. It outputs $False$ if and only if there exists an FTCG $\Gf$ in which there is a backdoor path from $X^i_{t -\gamma_i}$ to $Y_t$, otherwise it outputs $True$. Therefore, by Theorem \ref{th:equiv_IBC_multivarie}, Algorithm \ref{algo:calcul_IBC} is correct.

    Algorithm \ref{algo:calcul_IBC} runs Algorithm \ref{algo:calcul_t_NC} in $\mathcal{O}\left(\left| \Xf \right| \log \left| \Xf \right| + (\left| \mathcal{E}^s \right| + \left| \mathcal{V}^s \right|) \log \left| \mathcal{V}^s \right|\right)$, tests the existence of directed paths in $\mathcal{O}\left(\left| \mathcal{V}^s \right| + \left| \mathcal{E}^s \right| + |\Xf|\right)$, calls Algorithm \ref{algo:calcul_V_E_acc} $\left(\left| \Xf \right| + 1\right)$ times in $\mathcal{O}\left(\left| \Xf \right|(\left| \mathcal{E}^s \right| + \left| \mathcal{V}^s \right| \log \left| \mathcal{V}^s \right|)\right)$ and does $\mathcal{O}\left(|\Xf| \cdot \left| \mathcal{V}^s \right|\right)$ comparisons. Therefore, its overall complexity is $\mathcal{O}\left(\left| \Xf \right|\left( \log \left| \Xf \right|+ ( \left| \mathcal{E}^s \right| + \left| \mathcal{V}^s \right| )\log \left| \mathcal{V}^s \right|\right)\right)$.
\end{proof}

\newpage
\section{Proofs of Section \ref{sec:Consistency_Time}}\label{sec:proof:5}

\mylemmaequivexistencecheminforkNC*

\begin{figure}
    \centering
    \begin{subfigure}{0.3\textwidth}
        \centering
        \begin{tikzpicture}[scale = 1.5, ->,>=stealth,auto,node distance=3cm,semithick]
            \node at (0,0) (V1) {$V^1_{t_2}$};
            \node at (1,0) (V2) {$V^2_{t_2}$};
            \node at (1,1) (V3) {$V^2_{t_1}$};
            \node at (0,1) (V4) {$V^1_{t_1}$};
            \node at (0.5,2.1) (F) {$F_{t_{\NC}(F)}$};
            \node at (-0.5,-1) (X) {$X^i_{t- \gamma_i}$};
            \node at (1.5,-1) (Y) {$Y_t$};

            \path (V4) edge[CentraleBlue] (V3)
                  (F)  edge[CentraleBlue, decorate, decoration={snake, pre length = 3pt, post length=2pt, amplitude=1.5pt}] (V4)
                  (V3) edge[CentraleBlue, bend right = 40, decorate, decoration={snake, pre length = 3pt, post length=2pt, amplitude=1.5pt}] (X)
                  (V2) edge[CentraleRed]  (V1)
                  (F)  edge[CentraleRed, bend left = 45, decorate, decoration={snake, pre length = 3pt, post length=5pt, amplitude=1.5pt}] (V2)
                  (V1) edge[CentraleRed, decorate, decoration={snake, pre length = 3pt, post length=5pt, amplitude=1.5pt}]  (Y)
                  (V3) edge[olive] (V1);
        \end{tikzpicture}
        \caption{\label{subfig:CheminsForkNC:a}}
    \end{subfigure}
    \begin{subfigure}{0.3\textwidth}
        \centering
        \begin{tikzpicture}[scale = 1.5, ->,>=stealth,auto,node distance=3cm,semithick]
            \node at (0,1) (V1) {$V^1_{t_2}$};
            \node at (1,1) (V2) {$V^2_{t_2}$};
            \node at (1,0) (V3) {$V^2_{t_1}$};
            \node at (0,0) (V4) {$V^1_{t_1}$};
            \node at (0.5,2.1) (F) {$F_{t_{\NC}(F)}$};
            \node[outer sep = -5] at (-0.5,-1) (X) {$X^i_{t- \gamma_i}$};
            \node at (1.5,-1) (Y) {$Y_t$};

            \path (V4) edge[CentraleBlue] (V3)
                  (F)  edge[CentraleBlue, bend right = 45, decorate, decoration={snake, pre length = 3pt, post length=2pt, amplitude=1.5pt}] (V4)
                  (V3) edge[CentraleBlue, decorate, decoration={snake, pre length = 3pt, post length=2pt, amplitude=1.5pt}] (X)
                  (V2) edge[CentraleRed]  (V1)
                  (F)  edge[CentraleRed, decorate, decoration={snake, pre length = 3pt, post length=2pt, amplitude=1.5pt}] (V2)
                  (V1) edge[CentraleRed,bend left = 40, decorate, decoration={snake, pre length = 3pt, post length=5pt, amplitude=1.5pt}]  (Y)
                  (V1) edge[olive] (V3);
        \end{tikzpicture}
        \caption{\label{subfig:CheminsForkNC:b}}
    \end{subfigure}
    \newline
    
    \begin{subfigure}{0.3\textwidth}
    \centering
        \begin{tikzpicture}[scale = 2, ->,>=stealth,auto,node distance=3cm,semithick]
            \node[outer sep = -4] at (1,.7) (V1) {$V^1_{t-\gamma_i}$};
            \node at (1,-.2) (V2) {$V^1_{t_2}$};
            \node at (0,.7) (X1) {$X^i_{t-\gamma_i}$};
            \node at (0,-.2) (X2) {$X^i_{t_2}$};
            \node at (2,1) (F) {$F_{t_{\NC}(F)}$};
            \node at (2.5,-.5) (Y) {$Y_t$};

            \path (F)  edge[CentraleBlue, decorate, decoration={snake, pre length = 3pt, post length=5pt, amplitude=1.3pt}] (V1)
                  (V1) edge[CentraleBlue] (X1)
                  (F)  edge[CentraleRed, bend left = 18, decorate, decoration={snake, pre length = 3pt, post length=5pt, amplitude=1.3pt}] (X2)
                  (X2) edge[CentraleRed] (V2)
                  (V2) edge[CentraleRed, decorate, decoration={snake, pre length = 3pt, post length=5pt, amplitude=1.3pt}] (Y)
                  (V1)  edge[olive, bend right = 10, decorate, decoration={snake, pre length = 3pt, post length=5pt, amplitude=1.3pt}] (Y);

        \end{tikzpicture}
        \caption{\label{subfig:CheminsForkNC:c}}
    \end{subfigure}
    \begin{subfigure}{0.3\textwidth}
    \centering
        \begin{tikzpicture}[scale = 2, ->,>=stealth,auto,node distance=3cm,semithick]
            \node at (0,0) (V2) {$V^2_{t_2}$};
            \node at (1,0) (F2) {$F_{t_2}$};
            \node at (1,1) (F1) {$F_{t_{\NC}(F)}$};
            \node at (0,1) (V1) {$V^2_{t_{\NC}(F)}$};
            \node at (-.9,.7) (X) {$X^i_{t-\gamma_i}$};
            \node at (1.9,-.5) (Y) {$Y_t$};

            \path (F1) edge[CentraleBlue] (V1)
                  (V1) edge[CentraleBlue, decorate, decoration={snake, pre length = 3pt, post length=5pt, amplitude=1.3pt}] (X)
                  (F1) edge[CentraleRed, decorate, decoration={snake, pre length = 3pt, post length=5pt, amplitude=1.3pt}] (V2)
                  (V2) edge[CentraleRed] (F2)
                  (F2) edge[CentraleRed, decorate, decoration={snake, pre length = 3pt, post length=5pt, amplitude=1.3pt}] (Y)
                  (F1) edge [olive, decorate, decoration={snake, pre length = 3pt, post length=5pt, amplitude=1.3pt}] (Y);
        \end{tikzpicture}
        \caption{\label{subfig:CheminsForkNC:d}}
    \end{subfigure}
    \hfill
    \begin{subfigure}{0.3\textwidth}
    \centering
    \begin{tikzpicture}[scale = 2, ->,>=stealth,auto,node distance=3cm,semithick]
        \node at (1,1) (F) {$F_{t_{\NC}(F)}$};
        \node at (1,0) (V2) {$F_{t_2}$};
        \node at (0,1) (X1) {$X^i_{t-\gamma_i}$};
        \node at (0,0) (X2) {$X^i_{t_2}$};
        \node at (2,-.5) (Y) {$Y_t$};

        \path (F) edge[CentraleBlue] (X1)
              (F)  edge[CentraleRed,decorate, decoration={snake, pre length = 3pt, post length=5pt, amplitude=1.3pt}] (X2)
              (X2) edge[CentraleRed] (V2)
              (V2) edge[CentraleRed, decorate, decoration={snake, pre length = 3pt, post length=5pt, amplitude=1.3pt}] (Y)
              (F)  edge[olive, decorate, decoration={snake, pre length = 3pt, post length=5pt, amplitude=1.3pt}] (Y);
    \end{tikzpicture}
    \caption{\label{fig:CheminsForkNC:e}}
    \end{subfigure}
    
    \caption{ $\pi^f_1$ is represented in {\color{CentraleBlue} blue}, $\pi^f_2$ is represented in {\color{CentraleRed} red}, and the modification to be made to construct $\pi^f_3$ in $\Gf_3$ is represented in {\color{olive} green}.} \label{fig:CheminsForkNC}
\end{figure}
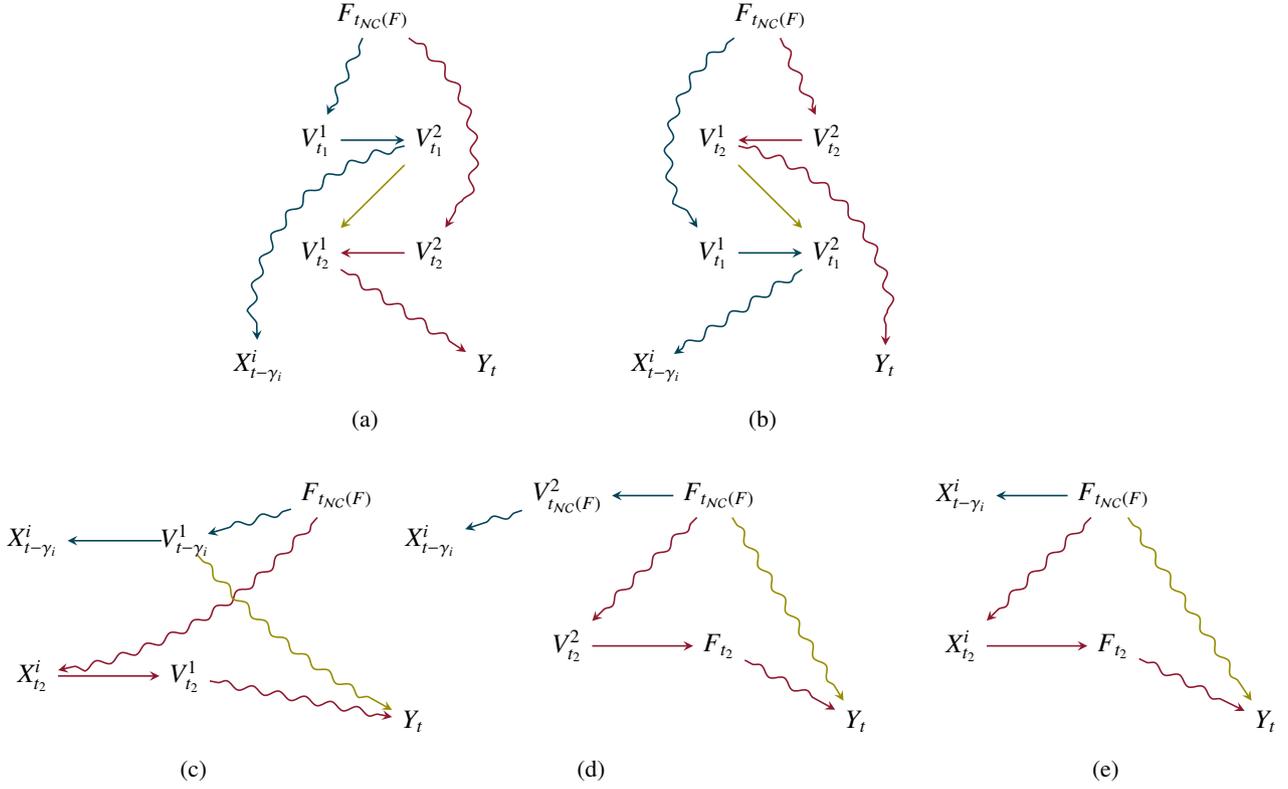

\begin{proof}
    Let us prove the direct implication (\ref{lemma:equiv_existence_chemin_fork_NC:1} $\Rightarrow$ \ref{lemma:equiv_existence_chemin_fork_NC:2}). Let $F$ and $X^i_{t - \gamma_i}$ be such that there exists an FTCG $\Gf \in \C(\Gs)$ that contains the path $X^i_{t- \gamma_i} \leftsquigarrow F_{t'} \rightsquigarrow Y_t$ which remains in $\NC \cup \{ X^i_{t-\gamma_i}\}$. Two cases arise:
    \begin{itemize}
        \item If $F \neq Y$ or $t - \gamma_i \neq t_{\NC}(F)$, the accessibility in $\Gf$ ensures that $F_{t_{\NC}(F)}$ is both $X^i_{t - \gamma_i}$-$\NC$-accessible and $Y_t$-$\NC$-accessible. Thus, we have proven Proposition \ref{lemma:equiv_existence_chemin_fork_NC:2a} in this case.
    
        \item Otherwise, $Y = F$ and $t - \gamma_i = t_{\NC}(F) = t_{\NC}(Y)$. Let us prove by contradiction that Lemma \ref{lemma:equiv_existence_chemin_fork_NC} \ref{lemma:equiv_existence_chemin_fork_NC:2b} holds. If it does not hold, then $\Gf$ must necessarily contain the edges $X^i_{t - \gamma_i} \leftarrow Y_{t_{\NC}(Y)}$ and $X^i_t \rightarrow Y_t$, which contradicts Assumption \ref{ass:Consistency_Time}.
    \end{itemize}

    \noindent Conversely, let us prove the indirect direction (\ref{lemma:equiv_existence_chemin_fork_NC:2} $\Rightarrow$ \ref{lemma:equiv_existence_chemin_fork_NC:1}). To prove this implication, we must show that  Lemma \ref{lemma:equiv_existence_chemin_fork_NC} \ref{lemma:equiv_existence_chemin_fork_NC:2a} implies  Lemma \ref{lemma:equiv_existence_chemin_fork_NC}  \ref{lemma:equiv_existence_chemin_fork_NC:1}, and that  Lemma \ref{lemma:equiv_existence_chemin_fork_NC}  \ref{lemma:equiv_existence_chemin_fork_NC:2b} implies  Lemma \ref{lemma:equiv_existence_chemin_fork_NC}  \ref{lemma:equiv_existence_chemin_fork_NC:1}. For each of these proofs, we are given $\Gf_1$, which contains $\pi^f_1 \vcentcolon= X^i_{t- \gamma_i} \leftsquigarrow F_{t_{\NC}(F)}$, which remains in $\NC \cup \{ X^i_{t- \gamma_i}\}$, and $\Gf_2$, which contains $\pi^f_2 \vcentcolon= F_{t_{\NC}(F)} \rightsquigarrow Y_t$, which remains in $\NC$. To prove Proposition \ref{lemma:equiv_existence_chemin_fork_NC:1}, it suffices to construct $\Gf_3$ that contains $\pi^f_3 \vcentcolon= X^i_{t- \gamma_i} \leftsquigarrow F_{t'} \rightsquigarrow Y_t$, which remains in $\NC \cup \{ X^i_{t- \gamma_i}\}$.

    \noindent Furthermore, without loss of generality, we can assume that the only intersection between $\pi^f_1$ and $\pi^f_2$ is $F_{t_{\NC}(F)}$. Indeed, consider $V_{t_v}$, the last element of $\pi^f_2$ in $\pi^f_1 \cap \pi^f_2$. By contradiction, we show that $V_{t_v} \neq Y_t$. If $V_{t_v} = Y_t$, then $\Gf_1$ contains $X^i_{t- \gamma_i} \leftsquigarrow Y_t$. However, by assumption, no FTCG contains a backdoor path without a fork from an intervention to $Y_t$. Thus, $V_{t_v} \neq Y_t$, and we can therefore work on ${\pi^f_1}' \vcentcolon= X^i_{t- \gamma_i} \leftsquigarrow V_{t_{\NC}(V)}$ and ${\pi^f_2}' \vcentcolon=V_{t_{\NC}(V)} \rightsquigarrow Y_t$. Therefore, concatenating $\pi^f_1$ and $\pi^f_2$ to create an FTCG does not create a cycle. Only Assumption \ref{ass:Consistency_Time} can be violated.

    If it is possible to concatenate $\pi^f_1$ and $\pi^f_2$ while maintaining Assumption \ref{ass:Consistency_Time}, then we prove Proposition \ref{lemma:equiv_existence_chemin_fork_NC:1}. If this is not the case, we show that, thanks to the assumptions of \ref{lemma:equiv_existence_chemin_fork_NC:2a} or \ref{lemma:equiv_existence_chemin_fork_NC:2b}, we can construct a suitable $\pi^f_3$. We thus merge the proofs of (\ref{lemma:equiv_existence_chemin_fork_NC:2a} $\Rightarrow$ \ref{lemma:equiv_existence_chemin_fork_NC:1}) and (\ref{lemma:equiv_existence_chemin_fork_NC:2b} $\Rightarrow$ \ref{lemma:equiv_existence_chemin_fork_NC:1}) as the reasoning is identical. 
    
    \noindent Therefore, we assume that it is not possible to concatenate $\pi^f_1$ and $\pi^f_2$ because Assumption \ref{ass:Consistency_Time} would be violated.
    We denote by $V^1_{t_1} \rightarrow V^2_{t_1}$ the last arrow of $\pi^f_1$ that contradicts $\pi^f_2$, and $V^1_{t_2} \leftarrow V^2_{t_2}$ the last arrow of $\pi^f_2$ that contradicts $V^1_{t_1} \rightarrow V^2_{t_1}$. We proceed by case distinction:
    \begin{itemize}
        \item If $V^1_{t_1} \neq F_{t_{\NC}(F)}$ and $V^2_{t_1} \neq X^i_{t- \gamma_i}$, we further distinguish three cases:
        \begin{itemize}
            \item If $t_1 < t_2$, then by adding the arrow $V^2_{t_1} \rightarrow V^1_{t_2}$, we can construct $\pi^f_3 = X^i_{t- \gamma_i} \leftsquigarrow V^2_{t_1}  \rightarrow V^1_{t_2} \rightsquigarrow Y_t$ (see Figure \ref{subfig:CheminsForkNC:a}).
            
            \item If $t_1 > t_2$, then by adding the arrow $V^1_{t_2} \rightarrow V^2_{t_1}$, we can construct $\pi^f_3 = X^i_{t- \gamma_i} \leftsquigarrow V^2_{t_1}  \leftarrow V^1_{t_2} \rightsquigarrow Y_t$ (see Figure \ref{subfig:CheminsForkNC:b}).

            \item The case $t_1 = t_2$ is excluded because the only intersection between $\pi^f_1$ and $\pi^f_2$ is $F_{t_{\NC}(F)}$.
        \end{itemize}
        
        \item If $V^1_{t_1} \neq F_{t_{\NC}(F)}$ and $V^2_{t_1} = X^i_{t- \gamma_i}$. $\Gf_2$ contains $V^1_{t_2} \rightsquigarrow Y_t$, which can be transformed into $V^1_{t-\gamma_i} \rightsquigarrow Y_t$ by changing the first arrow, yielding $\pi^f_3 = X^i_{t-\gamma_i} \leftarrow V^1_{t-\gamma_i} \rightsquigarrow Y_t$ (see Figure \ref{subfig:CheminsForkNC:c}).
        
        \item If $V^1_{t_1} = F_{t_{\NC}(F)}$ and $V^2_{t_1} \neq X^i_{t- \gamma_i}$. $\Gf_2$ contains $F_{t_2} \rightsquigarrow Y_t$, which can be transformed into $F_{t_{\NC}(F)}\rightsquigarrow Y_t$ by only changing the first arrow, yielding $\pi^f_3 = X^i_{t-\gamma_i} \leftsquigarrow V^2_{t_{\NC}(F)} \leftarrow F_{t_{\NC}(F)}  \rightsquigarrow Y_t$ (see Figure \ref{subfig:CheminsForkNC:d}).
        
        \item If $V^1_{t_1} = F_{t_{\NC}(F)}$ and $V^2_{t_1} = X^i_{t- \gamma_i}$, we further distinguish two cases:
        \begin{itemize}
            \item If $F \neq Y$, then $t - \gamma_i = t_{\NC(F)}$, because the case $t - \gamma_i \neq t_{\NC(F)}$ does not contradict Assumption \ref{ass:Consistency_Time}. $\Gf_2$ contains $F_{t_2} \rightsquigarrow Y_t$, which can be transformed into $F_{t_{\NC}(F)}\rightsquigarrow Y_t$, yielding $\pi^f_3 = X^i_{t-\gamma_i} \leftarrow V_{t_{\NC}(F)} \leftarrow F_{t_{\NC}(F)}  \rightsquigarrow Y_t$ (see Figure \ref{fig:CheminsForkNC:e}).
            
            \item Otherwise, if $F = Y$, the only case that can contradict Assumption \ref{ass:Consistency_Time} is $t - \gamma_i = t_{\NC(F)}$, i.e., we are strictly under the assumptions of Proposition \ref{lemma:equiv_existence_chemin_fork_NC:2b}. According to Proposition \ref{lemma:equiv_existence_chemin_fork_NC:2b}, $Y_{t_{\NC}(Y)}$ is $Y_t$ and is $\NC$-accessible without using $X^i_t \rightarrow Y_t$. Therefore, we can construct $\pi^f_3$ without contradicting Assumption \ref{ass:Consistency_Time}.
        \end{itemize}
    \end{itemize}

\end{proof}

\newpage
\section{Speed up of Algorithm \ref{algo:calcul_IBC}}\label{sct:speed_up_algo_3}

It is possible to compute 
$
\left( \max_{X^i_{t-\gamma_i}} \left\{t^{\NC}_{X^i_{t-\gamma_i}}(S) \right\}\right)_{S \in \mathcal{V}^s}
$
in a single traversal of the graph \(\Gs\). To achieve this, one only needs to modify the initialization of the priority queue \(Q\) in Algorithm \ref{algo:calcul_V_E_acc}. In this optimized strategy, \(Q\) is initialized with 
$
\left\{X^i_{t-\gamma_i}\right\}_{i \in \{1,\dots, n\}}
$. This is equivalent to computing 
$
\left(t^{\NC}_{M_{\star}}(S)\right)_{S \in \mathcal{V}^s},
$
where \(M_{\star}\) is a fictitious vertex corresponding to the merging of all interventions in on variable . The complexity remains in 
\(
\mathcal{O} \left(\left| \mathcal{E}^s \right| + \left| \mathcal{V}^s \right| \log \left| \mathcal{V}^s \right|\right),
\)
as the single traversal avoids redundant checks.

During this computation, it is also possible to efficiently test for the existence of collider-free backdoor paths involving a fork that remains within \(\NC\), except perhaps for its first vertex. Since 
\(
(t^{\NC}_{Y_{t}}(S))_{S \in \mathcal{V}^s}
\)
is already computed, the algorithm can, for each parent \(P\) encountered, directly verify whether 
\(
t_{\NC}(P) \leq t^{\NC}_{M_{\star}}(P) = \max_{X^i_{t-\gamma_i}} \left\{t^{\NC}_{X^i_{t-\gamma_i}}(P)\right\}
\)
and whether 
\(
t_{\NC}(P) \leq t^{\NC}_{Y_t}(F).
\)
This characterizes the existence of a collider-free backdoor path with a fork remaining in \(\CF\) from any intervention to \(Y_t\). Therefore, it is possible to test for the existence of a backdoor path with a fork in $\mathcal{O} \left(\left| \mathcal{E}^s \right| + \left| \mathcal{V}^s \right| \log \left| \mathcal{V}^s \right| \right)$.

Thus, the need for separate loops to test the existence of backdoor paths with forks is eliminated and identifiability by common adjustment can be computed in $\mathcal{O} \left(\left| \Xf \right| \log \left| \Xf \right| + \left| \mathcal{E}^s \right| + \left| \mathcal{V}^s \right| \log \left| \mathcal{V}^s \right| \right)$.

\section{An Efficient Algorithm for Identifiability by Common Adjustment Under Assumption \ref{ass:Consistency_Time}} \label{sct:IBC_pseudo_lineaire_consistent}

The algorithm presented in Algorithm \ref{algo:calcul_IBC_with_consistency} aims to determine whether the causal effect $P(y_t \mid \text{do}(x^1_{t-\gamma_1}), \dots, \text{do}(x^n_{t-\gamma_n}))$ is identifiable by common adjustment under Assumption \ref{ass:Consistency_Time}. Its correctness arises from Theorem \ref{th:equiv_IBC_multivarie}. Algorithm \ref{algo:calcul_IBC_with_consistency} tests the existence of an FTCG which contains a collider-free backdoor path from an intervention to $Y_t$ which remains in $\CF$. It starts by checking the existence of a directed path and then test the existence of backdoor paths with a fork using the characterisation of Lemma \ref{lemma:equiv_existence_chemin_fork_NC}. 

The algorithm checks the existence of a directed collider-free backdoor path that remains in $\CF$ as in Algorithm~\ref{algo:calcul_IBC}. The rest of the algorithm then focuses on checking the existence of backdoor path with a fork. To do so, it uses the characterisation of Lemma \ref{lemma:equiv_existence_chemin_fork_NC}. On line 1, it checks the backdoor path whose fork reduces to $F \neq Y$. Thanks to the strategy presented in Subsection \ref{sct:speed_up_algo_3}, this can be done in a single traversing. Then, on line 2, the algorithm checks that last case of case \ref{lemma:equiv_existence_chemin_fork_NC:2a}. It works on path whose fork reduces to $Y$ and whose interventions are at time $t- \gamma_i \neq t_{\NC}(Y)$. Thanks to the strategy presented in Subsection \ref{sct:speed_up_algo_3}, this can be done in a single traversing. Then, the algorithm checks the cases \ref{lemma:equiv_existence_chemin_fork_NC:2b}. It starts by the case \ref{lemma:equiv_existence_chemin_fork_NC:2b:i} on line 3. Since $\forall X^i_{t-\gamma_i} \in \mathcal{X}, t-\gamma_i = t_{\NC}(Y)$, finding a directed path from $Y_{t_{\NC}(Y)}$ to $X^i_{t_{\NC}(Y)}$ without using $X^i_{t_{\NC}(Y)} \leftarrow Y_{t_{\NC}(Y)} $ is equivalent of finding a directed path from $Y$ to $X^i$ in $\Gs$ without using $Y \rightarrow X^i$. This can be done by a BFS algorithm in $\Gs$. The algorithm finishes by checking the case \ref{lemma:equiv_existence_chemin_fork_NC:2b:ii}. $\mathcal{X}'$ represents the set of interventions at time $t_{\NC}(Y)$ for which $Y_{t_{\NC}(Y)}$ is $X^i_{t-\gamma_i}$-$\NC$-accessible. If $\left| \mathcal{X}' \right| \geq 2$, since $Y_{t_{\NC}(Y)}$ is $Y_t$-$\NC$-accessible, there is at least one $X^i_{t-\gamma_i} \in \mathcal{X}'$ that $Y_{t_{\NC}(Y)}$ is $Y_t$ - $\NC$-accessible without using $X^i_t \rightarrow Y_t$. Indeed,  $Y_{t_{\NC}(Y)}$ is $Y_t$-$\NC$-accessible, therefore, there exists an FTCG which contains a directed path from $Y_{t_{\NC}(Y)}$ to $Y_t$. This path cannot use simultaneously an arrow $X^i_t \rightarrow Y_t$ and an arrow $X^j_t \rightarrow Y_t$.  This reasoning explains the fourth line of the algorithm. The line 5 is reached when $\left| \mathcal{X}' \right| = 1$. In this case, a traversing algorithm starting from $Y_{t_{\NC}(Y)}$ can test if $Y_{t_{\NC}(Y)}$ is $Y_t$ - $\NC$-  accessible without using $X^i_t \rightarrow Y_t$. All the subcases of case  \ref{lemma:equiv_existence_chemin_fork_NC:2b:ii} have been checked by the algorithm

Algorithm \ref{algo:calcul_IBC_with_consistency} calls 5 traversing of $\Gs$. Therefore, its complexity is indeed pseudo-linear.

\begin{algorithm}[t]

\caption{Computation of identifiability by common adjustment under Assumption \ref{ass:Consistency_Time}.}
\label{algo:calcul_IBC_with_consistency}
\SetNlSty{}{\relsize{2.5}}{.}
\SetKwInput{KwData}{Input}
\SetKwInput{KwResult}{Output}

\KwData{ $\Gs$ an SCG and $P(y_t \mid \text{do}(x^1_{t-\gamma_1}), \dots, \text{do}(x^n_{t-\gamma_n}))$ the considered effect.}
\KwResult{A boolean indicating whether the effect is identifiable by common adjustment or not.}

$(t_{\NC}(S))_{S \in \mathcal{V}^s} \gets $ Algorithm \ref{algo:calcul_t_NC} \;
\tcp{Enumeration of directed paths.} 
$\mathcal{S} \gets \{S \in \mathcal{V}^s \mid t_{\NC}(S) \leq 0 \}\cup \{X^i \mid t-\gamma_i = 0\}$ \;
\If{$\exists i \in \{1,\ldots, n\}$ s.t. $ X^i \in \Desc \left(Y, \Gs_{\mid \mathcal{S}} \right)$ and $\gamma_i = 0$}{
    \Return{False}
}

\tcp{Enumeration of fork paths.} 

$(t^{\NC}_{Y_{t}}(S))_{S \in \mathcal{V}^s} \gets $ Algorithm \ref{algo:calcul_V_E_acc} \;

\tcp{Test all all forks $F \neq Y$:}

\nl \If(\tcc*[f]{\ref{sct:speed_up_algo_3}}){there exist $X^i_{t- \gamma_i}$ and $F$ such that $F_{t_{\NC}(F)}$ is well defined, $X^i_{t- \gamma_i}$-$\NC$-accessible and $Y_t$-$\NC$-accessible}{
    \Return{False}
}

\tcp{Test for $F = Y$:}
\If{$Y_{t_{\NC}(Y)}$ is $Y_t$-$\NC$-accessible}{
    \tcp{Last case of \ref{lemma:equiv_existence_chemin_fork_NC:2a}}
    
    \nl \If(\tcc*[f]{\ref{sct:speed_up_algo_3}}){there exist $X^i_{t- \gamma_i}$ such that $t- \gamma_i \neq t_{\NC}(Y)$ and $Y_{t_{\NC}(Y)}$ is $X^i_{t- \gamma_i}$-$\NC$-accessible}{ 
        \Return{False}
    }
    
    \tcp{Test for \ref{lemma:equiv_existence_chemin_fork_NC:2b}:}
    
    $\mathcal{X} \gets \{X^i_{t-\gamma_i} \mid t-\gamma_i = t_{\NC}(Y) \}$ \;
    
    \tcp{Case \ref{lemma:equiv_existence_chemin_fork_NC:2b:i}:}
    
    \nl \If(\tcc*[f]{BFS}){ $\exists~ X^i_{t-\gamma_i} \in \mathcal{X}$ such that $Y_{t_{\NC}(Y)}$ is $X^i_{t- \gamma_i}$ - $\NC$-accessible without using $X^i_{t- \gamma_i} \leftarrow Y_{t- \gamma_i}$ } {
        \Return{False}
    }
    \tcp{Case \ref{lemma:equiv_existence_chemin_fork_NC:2b:ii}:}
    $\mathcal{X}' \gets \{X^i_{t-\gamma_i} \mid t-\gamma_i = t_{\NC}(Y) \text{ and } X^i \in \Desc(Y, \Gs) \}$ \tcc*{BFS}

    \nl \If{$\left| \mathcal{X}' \right| \geq 2$} {
        \Return{False}
    } 
    \If(\tcc*[f]{BFS}){$\mathcal{X}' = \{ X^i_{t-\gamma_i}\}$ \textbf{and} $Y_{t_{\NC}(Y)}$ is $Y_t$ - $\NC$-  accessible without using $X^i_t \rightarrow Y_t$.} {
        \Return{False}
        
    }
}
\Return{True}
\end{algorithm}

\end{document}